\documentclass{amsart}
\usepackage{amsfonts,amsmath,amssymb,mathrsfs,array}
\usepackage[english]{babel}
\usepackage[hidelinks]{hyperref}
\newtheorem{thm}{Theorem}[section]
\newtheorem{lem}[thm]{Lemma}

\newtheorem{proposition}[thm]{Proposition}

\newtheorem{corollary}[thm]{Corollary}
\theoremstyle{definition}
\newtheorem{definition}[thm]{Definition}
\newtheorem{rem}[thm]{Remark}
\newtheorem{remark}[thm]{Remark}
\numberwithin{equation}{section}

\newcommand{\1}{\mathbf{1}}

\begin{document}
\title[Convolution on quantum tori]
{Convolution on quantum tori I: sharp Young inequalities and endpoint Sobolev embeddings}
\author[F. Sukochev]{F. Sukochev}
\address{School of Mathematics and Statistics, University of New South Wales,
Kensington, NSW 2052, Australia}
\email{f.sukochev@unsw.edu.au}

\author[R. Tastankul]{R. Tastankul}
\address{Institute of Mathematics and Mathematical Modeling,
050010 Almaty, Kazakhstan}
\email{ramazan.tastankul@gmail.com}

\author[K. Tulenov]{K. Tulenov}
\address{ School of Mathematics and Statistics,
	University of New South Wales, Kensington, NSW 2052, Australia; Institute of Mathematics and Mathematical Modeling,
050010 Almaty, Kazakhstan}
\email{kanat.tulenov@unsw.edu.au}

\author[D. Zanin]{D. Zanin}
\address{School of Mathematics and Statistics, Central South University,
Changsha 410075, China}
\email{d.zanin@csu.edu.cn}

\date{}

\begin{abstract}
We introduce a convolution of  noncommutative tori.  The construction is the preadjoint of a trace-preserving cocycle-addition homomorphism. Its zero-parameter case is the function--operator semi-convolution. We prove both mixed $L_{1}$--$L_\infty$ inequalities and various Young inequalities with the optimal constant one, including their limiting cases. The endpoint bounds also yield submajorization and convolution estimates in fully symmetric spaces. As applications, we establish, for every integer $k\geq1$ and $1\leq p\leq\infty$, $W_{p}^{k}(\mathbb{T}_{\theta}^{d})\hookrightarrow B_{\infty,\infty}^{k-\frac{d}{p}}(\mathbb{T}_{\theta}^{d})$, and prove the critical embedding $W_{1}^{d}(\mathbb{T}_{\theta}^{d})\hookrightarrow L_\infty(\mathbb{T}_{\theta}^{d})$. These results settle the two questions posed in Remark~6.8 of Xiong, Xu and Yin \cite{XiongXuYin}. 
\end{abstract}



\maketitle
\tableofcontents

\section{Introduction}

The noncommutative, or quantum, torus is a fundamental example in Connes' noncommutative geometry \cite{Connes1980,ConnesNCGeometry}. Its structure as a noncommutative differentiable manifold was studied by Rieffel \cite{RieffelTori}. Quantum tori also occur in mathematical physics: Connes, Douglas and Schwarz \cite{ConnesDouglasSchwarz} used them in toroidal compactifications of Matrix theory, while Landi, Lizzi and Szabo \cite{LandiLizziSzabo} studied their role in string geometry. These geometric and physical developments have been accompanied by the study of integration, differentiation and Fourier analysis on quantum tori.  

For a real skew-symmetric $d\times d$ matrix $\theta$, the $C^*$-algebra $C(\mathbb{T}_\theta^d)$ is generated by unitaries $U_1,\ldots,U_d$ satisfying
\begin{equation*}
U_jU_k=e^{2\pi i\theta_{jk}}U_kU_j,
\ \  1\leq j,k\leq d.
\end{equation*}
When $\theta=0$, this is the algebra $C(\mathbb{T}^d)$ of continuous functions on the ordinary torus. The canonical trace $\tau_\theta$ gives a tracial von Neumann algebra $L_\infty(\mathbb{T}_\theta^d)$ and the associated noncommutative spaces $L_p(\mathbb{T}_\theta^d)$, with $\tau_\theta(\mathbf1)=1$. For the construction and differential calculus, we refer to \cite[Section~3.5]{LordMcDonaldSukochevZanin}. Spera \cite{Spera} studied Sobolev spaces on noncommutative tori. Chen, Xu and Yin \cite{ChenXuYin} developed Fourier summation, Fourier multipliers and Littlewood--Paley theory in this setting; Xiong, Xu and Yin \cite{XiongXuYin} established a systematic theory of Sobolev, Besov and Triebel--Lizorkin spaces. Further results on multipliers, quantum differentiability and applications to evolution equations appear in \cite{Lai, McDonaldSukochevXiong, Ricard, RuzhanskyShaimardanTulenovTori, ShaimardanTastankulTulenov, SukochevTulenovZanin}. 

The purpose of the present paper is to introduce a convolution of two operators on quantum tori and to establish its Young inequalities and endpoint Sobolev applications. Recall that, on the classical torus $\mathbb{T}^d$ with normalized Haar measure, Young's inequality \cite{Young} states that
\begin{equation}\label{eq:intro-Young-relation}
\begin{gathered}
\|f\ast g\|_{L_r(\mathbb{T}^d)}
\leq\|f\|_{L_p(\mathbb{T}^d)}\|g\|_{L_q(\mathbb{T}^d)},\\
1+\frac1r=\frac1p+\frac1q,
\ \  1\leq p,q,r\leq\infty.
\end{gathered}
\end{equation}
On a quantum torus, a classical function $f$ can be convolved with an operator $x$ by averaging the canonical torus action:
\begin{equation*}
f\ast x=\int\limits_{\mathbb{T}^d}f(s)\alpha_{-s}^\theta(x)\,dm(s),
\ \  f\in L_1(\mathbb{T}^d),\ \  x\in L_1(\mathbb{T}_\theta^d).
\end{equation*}
Here $\mathbb{T}^d=\mathbb{R}^d/(2\pi\mathbb{Z})^d$, $m$ is normalized Haar measure, and $\alpha_s^\theta$ is determined by $\alpha_s^\theta(U_j)=e^{is_j}U_j$; see \cite[Section~2.3]{ChenXuYin} and \cite[Section~4]{STTZ}. We call this function--operator operation \emph{semi-convolution}. Our previous work \cite{STTZ} established an O'Neil-type rearrangement inequality for semi-convolution and applied it to fractional Sobolev embeddings in noncommutative symmetric spaces. The present paper directly continues that work.

We now allow both inputs to be operators, possibly on different quantum tori. The guiding identity is the classical Fourier formula $\widehat{f\ast g}(n)=\widehat f(n)\widehat g(n)$. For real skew-symmetric matrices $\theta_1$ and $\theta_2$, we construct a bilinear map
\begin{equation}\label{eq:intro-convolution}
L_1(\mathbb{T}_{\theta_1}^d)\times L_1(\mathbb{T}_{\theta_2}^d)
\to L_1(\mathbb{T}_{\theta_1+\theta_2}^d), \ \
(x,y)\mapsto x\star_{\theta_1,\theta_2}y,
\end{equation}
characterized by
\begin{equation}\label{eq:intro-fourier-convolution}
\widehat{x\star_{\theta_1,\theta_2}y}(n)
=\widehat x(n)\widehat y(n),\ \  n\in\mathbb{Z}^d.
\end{equation}
The Fourier coefficients here are taken with respect to the symmetrized Fourier generators specified in Section~\ref{sec:preliminaries}. For $\theta_1=0$, this convolution reduces to semi-convolution; when both parameters vanish, it is ordinary convolution on $\mathbb{T}^d$.

The addition of the parameters follows directly from the defining
commutation relations. Indeed, the unitaries
$U_j^{\theta_1}\otimes U_j^{\theta_2}$ satisfy those relations
with parameter $\theta_1+\theta_2$. They therefore determine a
trace-preserving normal $*$-homomorphism
\begin{equation*}
\mathcal T_{\theta_1,\theta_2}:
L_\infty(\mathbb{T}_{\theta_1+\theta_2}^d)
\rightarrow
L_\infty(\mathbb{T}_{\theta_1}^d)\,\overline\otimes\,
L_\infty(\mathbb{T}_{\theta_2}^d).
\end{equation*}
We define $x\star_{\theta_1,\theta_2}y$ by applying its
preadjoint to $x\otimes y$. A corresponding homomorphism was
used by Ricard \cite[proof of Proposition~2.1]{Ricard} in the
study of Fourier multipliers. Related tensor-product maps occur
in quantum Euclidean harmonic analysis: the corepresentation
and diagonal homomorphisms of Gonz\'alez-P\'erez, Junge and Parcet
\cite[Corollary~1.4 and Lemma~1.5]{GonzalezPerezJungeParcet}
connect a quantum Euclidean algebra with a classical factor or
with its opposite algebra. In the present construction, opposite
parameters $\theta_2=-\theta_1$ give an output on the classical
torus.

The construction immediately yields positivity and the
$L_1\times L_1\to L_1$ contraction. The main analytic step is
to prove the two mixed endpoint bounds
\begin{equation}\label{eq:intro-mixed}
\begin{split}
\|x\star_{\theta_1,\theta_2}y\|_\infty
&\leq\|x\|_1\|y\|_\infty,\\
\|x\star_{\theta_1,\theta_2}y\|_\infty
&\leq\|x\|_\infty\|y\|_1.
\end{split}
\end{equation}
Together with the $L_1\times L_1\to L_1$ contraction, these estimates give
Young type inequality, Theorem~\ref{thm:Young}:
\begin{equation}\label{eq:intro-sharp-Young}
\|x\star_{\theta_1,\theta_2}y\|_{L_r(\mathbb{T}_{\theta_1+\theta_2}^d)}
\leq\|x\|_{L_p(\mathbb{T}_{\theta_1}^d)}
     \|y\|_{L_q(\mathbb{T}_{\theta_2}^d)},
\end{equation}
whenever the exponents satisfy \eqref{eq:intro-Young-relation}.
The constant one is optimal, as is seen by taking
$x=y=\mathbf1$. The endpoint estimates also imply
\begin{equation}\label{eq:intro-HLP}
\begin{split}
\mu(x\star_{\theta_1,\theta_2}y)
&\prec\!\prec\|x\|_1\mu(y),\\
\mu(x\star_{\theta_1,\theta_2}y)
&\prec\!\prec\|y\|_1\mu(x),
\end{split}
\end{equation}
where $\mu$ denotes generalized singular values and
$\prec\!\prec$ is Hardy--Littlewood--P\'olya submajorization.
Thus the limiting convolution estimates extend to fully
symmetric operator spaces, using the noncommutative integration
framework of \cite{FackKosaki,DoddsDePagterSukochev}.
All these estimates are uniform in the deformation parameters.

Using these convolution estimates and explicit Fourier
reconstructions, we give affirmative answers to the two
questions posed by Xiong, Xu and Yin in
\cite[Remark~6.8]{XiongXuYin}. The first asks whether the
classical critical embedding remains valid on a quantum torus:
\begin{equation}\label{eq:intro-critical-Sobolev}
W_1^d(\mathbb{T}_\theta^d)
\hookrightarrow L_\infty(\mathbb{T}_\theta^d).
\end{equation}
The second concerns the Sobolev--Besov embedding at $p=1$
for odd integer orders. The derivative-Sobolev conclusion of
\cite[Theorem~6.6(ii)]{XiongXuYin} covers $1<p<\infty$ and,
when $p=1$, even orders. We prove
\begin{equation}\label{eq:intro-Sobolev-Besov}
W_p^k(\mathbb{T}_\theta^d)
\hookrightarrow B_{\infty,\infty}^{k-\frac{d}{p}}(\mathbb{T}_\theta^d),
\ \  k\geq1,\ \  1\leq p\leq\infty,
\end{equation}
for every integer $k$, including odd $k$ when $p=1$; see
Theorems~\ref{thm:W1-Besov-endpoint}
and~\ref{thm:Wp-Besov-odd-endpoint}.
The critical inclusion \eqref{eq:intro-critical-Sobolev} is
proved in Theorem~\ref{thm:critical-W1-Linfty}.
The spaces $W_p^k$ are defined by the canonical derivatives,
and $B_{\infty,\infty}^s$ by a Littlewood--Paley decomposition;
these conventions are fixed in the preliminaries.

For \eqref{eq:intro-Sobolev-Besov}, we reconstruct the nonzero
Fourier modes from the pure derivatives of order $k$ using
scalar kernels and estimate their dyadic blocks by convolution.
For \eqref{eq:intro-critical-Sobolev}, modes with all coordinates nonzero are reconstructed from the mixed derivative using a bounded periodic sawtooth kernel, whereas modes with vanishing coordinates are reconstructed via a finite decomposition. The resulting embedding constants
are independent of $\theta$.

The paper is organized as follows. Section~\ref{sec:preliminaries}
contains the preliminaries on quantum tori and function spaces.
Section~\ref{sec:convolution} constructs convolution, establishes
its algebraic properties and proves the mixed endpoint bounds.
Section~\ref{sec:young} proves Young's inequality, and
Section~\ref{sec:limiting} gives the limiting and fully symmetric
consequences. Section~\ref{sec:endpoint-sobolev} treats the two
Sobolev questions. The anti-isomorphism used in the endpoint
argument is constructed in Appendix~\ref{app:parameter-reversal}.

We are grateful to Professor Xiao Xiong for suggesting the two-operator construction and drawing our attention to its physical motivation.

\section{Preliminaries}\label{sec:preliminaries}

Throughout, $d\geq1$. We use the classical torus $\mathbb{T}^d=\mathbb{R}^d/(2\pi\mathbb{Z})^d$ with normalized Haar measure $m$; consequently $\widehat f(n)=\int\limits_{\mathbb{T}^d}f(t)e^{-in\cdot t}\,dm(t)$. A deformation parameter is a real skew-symmetric $d\times d$ matrix $\theta$. When $d=1$, necessarily $\theta=0$.

\subsection{Quantum tori and Fourier systems}

The symbol $\mathbb{T}_\theta^d$ denotes the quantum torus;
we distinguish its $C^*$-algebra from its von Neumann
algebra. Let $C(\mathbb{T}_\theta^d)$ be the universal
$C^*$-algebra generated by unitaries
$U_1,\ldots,U_d$ satisfying
\begin{equation}\label{eq:prelim-generators}
U_jU_k=e^{2\pi i\theta_{jk}}U_kU_j,
\ \  1\leq j,k\leq d.
\end{equation}
The ordered monomials $U^n=U_1^{n_1}\cdots U_d^{n_d}$,
$n\in\mathbb{Z}^d$, span a norm-dense $*$-subalgebra.
The prescription
$\tau_\theta(U^n)=\delta_{n,0}$ extends to a faithful
tracial state on $C(\mathbb{T}_\theta^d)$. In its GNS
representation, the weak operator closure is the finite
von Neumann algebra
\begin{equation}\label{eq:prelim-vna}
L_\infty(\mathbb{T}_\theta^d)
:=\pi_{\tau_\theta}\bigl(C(\mathbb{T}_\theta^d)\bigr)'',
\ \  \tau_\theta(\mathbf1)=1.
\end{equation}
The same symbol $\tau_\theta$ denotes the \emph{normal}
faithful tracial state on this von Neumann algebra.
For $\theta=0$, the two algebras are
$C(\mathbb{T}^d)$ and $L_\infty(\mathbb{T}^d)$,
respectively. This separation of the two algebras is
used throughout. See
\cite[Section~1.2]{XiongXuYin},
\cite[Section~2]{Ricard}, and
\cite[Section~3.5.1]{LordMcDonaldSukochevZanin};
the deformation parameter in the last reference uses
a different $2\pi$ normalization.

Our symmetrized Fourier generators are
\begin{equation}\label{eq:prelim-Weyl-definition}
e_n^\theta
=\exp\!\left(-\pi i\sum_{1\leq j<k\leq d}
\theta_{jk}n_jn_k\right)U^n,
\ \  n\in\mathbb{Z}^d.
\end{equation}
Writing $\langle n,\theta m\rangle
=\sum_{j,k}n_j\theta_{jk}m_k$, one obtains
\begin{equation}\label{eq:prelim-Weyl}
e_n^\theta e_m^\theta
=\sigma_\theta(n,m)e_{n+m}^\theta,
\ \ 
\sigma_\theta(n,m)
=e^{\pi i\langle n,\theta m\rangle},
\ \ 
(e_n^\theta)^*=e_{-n}^\theta.
\end{equation}
In particular, $\tau_\theta(e_n^\theta)=\delta_{n,0}$.
The finite span $\mathcal P_\theta$ of the
$e_n^\theta$ is the algebra of Fourier polynomials.
Its closure in $L_2(\mathbb{T}_\theta^d)$ has
$\{e_n^\theta:n\in\mathbb{Z}^d\}$ as an orthonormal
basis. For $x\in L_1(\mathbb{T}_\theta^d)$, put
\begin{equation}\label{eq:prelim-Fourier}
\widehat x(n)=\tau_\theta(xe_{-n}^\theta),
\ \ 
\|x\|_2^2=\sum_{n\in\mathbb{Z}^d}
|\widehat x(n)|^2\ \ (x\in L_2).
\end{equation}
Fourier coefficients determine elements of $L_1$
uniquely. Fourier polynomials are norm dense in
$L_p$, $1\leq p<\infty$, and weak$^*$ dense
in $L_\infty$; see \cite[Section~1]{ChenXuYin}
and \cite[Section~1]{Ricard}. The symmetrization in
\eqref{eq:prelim-Weyl-definition} ensures
$\sigma_{\theta_1+\theta_2}
=\sigma_{\theta_1}\sigma_{\theta_2}$.

\subsection{Measurable operators and symmetric spaces}

Let $(\mathcal M,\tau)$ be a finite von Neumann algebra
with a faithful normal tracial state, represented on a
Hilbert space $H$. An operator $x$ is affiliated with
$\mathcal M$ if it is closed, densely defined, and
commutes with every unitary in the commutant
$\mathcal M'$ (with the usual domain convention).
It is $\tau$-measurable if, for each $\varepsilon>0$,
there is a projection $e\in\mathcal M$ such that
$eH\subset\operatorname{Dom}(x)$ and
$\tau(\mathbf1-e)<\varepsilon$.
Write $S(\mathcal M,\tau)$ for these operators.
For the torus we abbreviate
$S(\mathbb{T}_\theta^d)
:=S(L_\infty(\mathbb{T}_\theta^d),\tau_\theta)$.
These operators need not belong to
$L_\infty(\mathbb{T}_\theta^d)$; affiliation,
measurability, and boundedness will not be identified.
See \cite{FackKosaki},
\cite[Chapter~1]{DoddsDePagterSukochev},
and \cite[Section~2.3]{STTZ}.

For $x\in S(\mathcal M,\tau)$ define its distribution
function and generalized singular values by
\begin{equation}\label{eq:prelim-mu}
d_x(\lambda)
=\tau\bigl(\mathbf1_{(\lambda,\infty)}(|x|)\bigr),
\ \ 
\mu(t,x)
=\inf\{\lambda\geq0:d_x(\lambda)\leq t\},
\ \  \lambda\geq0,\ 0<t<1.
\end{equation}
Here, $|x|=(x^\ast x)^{\frac{1}{2}}$ and
$\mu(\cdot,x)$ is decreasing. For a scalar function
$f$, the notation $\mu(\cdot,f)$ means its decreasing
rearrangement. The noncommutative Lebesgue spaces are
\begin{equation}\label{eq:prelim-Lp}
L_p(\mathcal M,\tau)
=\{x\in S(\mathcal M,\tau):\tau(|x|^p)<\infty\},
\ \ 
\|x\|_p=\tau(|x|^p)^{1/p}
\ \ (1\leq p<\infty),
\end{equation}
and $L_\infty(\mathcal M,\tau)=\mathcal M$ with the
operator norm. We write $L_p(\mathbb{T}_\theta^d)$
when $\mathcal M=L_\infty(\mathbb{T}_\theta^d)$.
The normalized trace gives
$L_\infty(\mathbb{T}_\theta^d)
\subset L_p(\mathbb{T}_\theta^d)
\subset L_1(\mathbb{T}_\theta^d)$.

We use Hardy--Littlewood--P\'olya submajorization:
\begin{equation}\label{eq:prelim-HLP}
x\prec\!\prec y
\ \ \Longleftrightarrow\ \ 
\int\limits_0^t\mu(s,x)\,ds
\leq\int\limits_0^t\mu(s,y)\,ds
\ \ (0<t\leq1).
\end{equation}
If $E$ is a symmetric Banach function space on
$(0,1)$, define its associated operator space by
\begin{equation}\label{eq:prelim-symmetric-space}
\mathcal E(\mathcal M,\tau)
=\{x\in S(\mathcal M,\tau):
      \mu(\cdot,x)\in E\},
\ \ 
\|x\|_{\mathcal E}
=\|\mu(\cdot,x)\|_E.
\end{equation}
As with $L_p$, we abbreviate
$\mathcal E(L_\infty(\mathbb{T}_\theta^d),\tau_\theta)$
to $\mathcal E(\mathbb{T}_\theta^d)$.
The space $E$ is fully symmetric if
$g\prec\!\prec f$ and $f\in E$ imply
$g\in E$ and $\|g\|_E\leq\|f\|_E$.
On a probability space,
$L_\infty\subset E\subset L_1$ continuously.
We will use the exact $K$-functional formula
\begin{equation}\label{eq:K-prelim}
K\bigl(t,x;L_1(\mathcal M,\tau),
             L_\infty(\mathcal M,\tau)\bigr)
=\int\limits_0^t\mu(s,x)\,ds,
\ \  0<t\leq1.
\end{equation}
These facts and their operator-space versions can be
found in \cite{FackKosaki,DoddsDePagterSukochev}.

\subsection{Distributions, Fourier multipliers, and derivations}

The smooth quantum torus is the Fr\'echet $*$-algebra
\begin{equation}\label{eq:prelim-smooth}
\mathcal S(\mathbb{T}_\theta^d)
=\left\{
\sum_{n\in\mathbb{Z}^d}a_ne_n^\theta:
\sup_{n\in\mathbb{Z}^d}(1+|n|)^N|a_n|<\infty
\text{ for every }N\geq0
\right\}.
\end{equation}
Its continuous dual
$\mathcal S'(\mathbb{T}_\theta^d)$
consists precisely of formal Fourier series with
coefficients of at most polynomial growth. With
the pairing $F_x(a)=\tau_\theta(xa)$ for
$x\in L_1$ and $a\in\mathcal S$, the
distributional Fourier coefficient is
$\widehat F(n)=F(e_{-n}^\theta)$.
Consequently
$L_p(\mathbb{T}_\theta^d)\subset
\mathcal S'(\mathbb{T}_\theta^d)$ for
$1\leq p\leq\infty$.
See \cite[Section~2]{XiongXuYin} and
\cite[Section~2.4]{STTZ}.

A scalar symbol $m:\mathbb{Z}^d\to\mathbb C$ of at
most polynomial growth determines the Fourier
multiplier
\begin{equation}\label{eq:prelim-multiplier}
T_m^\theta F
=\sum_{n\in\mathbb{Z}^d}
 m(n)\widehat F(n)e_n^\theta,
\ \  F\in\mathcal S'(\mathbb{T}_\theta^d).
\end{equation}
This formula defines a continuous map on
$\mathcal S'$ and restricts to a continuous map on
$\mathcal S$. When its restriction to
$\mathcal P_\theta$ extends boundedly to
$L_p(\mathbb{T}_\theta^d)$, we call $m$ an
$L_p$ Fourier multiplier. For $p=\infty$
the extension is understood to be normal and
weak$^*$ continuous; for $1\leq p<\infty$
it is uniquely determined by density. A bounded
symbol always gives an $L_2$ multiplier with
$\|T_m^\theta\|_{2\to2}=\|m\|_{\ell_\infty}$,
but boundedness of $m$ alone does not assert
$L_p$ boundedness at other exponents.
See \cite[Section~2.5]{STTZ},
\cite[Section~1]{Ricard}, and
\cite[Section~3.5.2]{LordMcDonaldSukochevZanin}.
The completely bounded multiplier norm, when
applicable, is independent of $\theta$
\cite[Proposition~2.1]{Ricard}; this is a
statement about \emph{complete} boundedness.

The canonical derivations and Laplacian are
the distributional Fourier multipliers
\begin{equation}\label{eq:canonical-derivations}
\partial_j^\theta e_n^\theta=in_je_n^\theta,
\ \ 
\Delta_\theta e_n^\theta=-|n|^2e_n^\theta,
\ \ 
J_\theta^s e_n^\theta
=(1+|n|^2)^{s/2}e_n^\theta.
\end{equation}
Thus
$\Delta_\theta=\sum_{j=1}^d(\partial_j^\theta)^2$
and $D_j:=-i\partial_j^\theta$ has symbol
$n_j$. The Fourier multiplier formulas define
all three operators on $\mathcal S'$; on $L_2$,
$-\Delta_\theta$ is positive and self-adjoint
with eigenvalues $|n|^2$.
Our angular normalization produces the
factor $i n_j$ rather than $2\pi i n_j$;
compare \cite[Section~2]{XiongXuYin},
\cite[Section~2.6]{STTZ}, and
\cite[Section~3.5.2]{LordMcDonaldSukochevZanin}.

The trace-preserving torus action
$\alpha_s^\theta$ is determined by
\begin{equation}\label{eq:prelim-action}
\alpha_s^\theta(e_n^\theta)
=e^{in\cdot s}e_n^\theta,
\ \  s\in\mathbb{T}^d.
\end{equation}
The action extends isometrically to $L_p$.
For $1\leq p<\infty$, each orbit
$s\mapsto\alpha_s^\theta(x)$ is
$L_p$-norm continuous; on $L_\infty$
the action is weak$^*$ continuous.
For $f\in L_1(\mathbb{T}^d)$ and
$x\in L_1(\mathbb{T}_\theta^d)$ define the
function--operator \emph{semi-convolution} by
\begin{equation}\label{eq:prelim-semi-convolution}
f\ast x
=\int\limits_{\mathbb{T}^d}
 f(s)\alpha_{-s}^\theta(x)\,dm(s),
\ \ 
\widehat{f\ast x}(n)
=\widehat f(n)\widehat x(n).
\end{equation}
The integral is a Bochner integral in $L_1$.
It is also defined in $L_p$ for
$1\leq p<\infty$ and $x\in L_p$,
and as a weak$^*$ integral for bounded $x$.
In particular, $f\ast x=T_{\widehat f}^\theta x$.
The minus sign in \eqref{eq:prelim-semi-convolution}
follows from our Fourier convention; compare
\cite[Section~2.5]{STTZ} and
\cite[Section~2.2]{McDonaldSukochevXiong}.

\subsection{Sobolev and Besov spaces}

For an integer $k\geq0$ and $1\leq p\leq\infty$,
the derivative Sobolev space is
\begin{equation}\label{eq:prelim-derivative-space}
W_p^k(\mathbb{T}_\theta^d)
=\left\{
x\in\mathcal S'(\mathbb{T}_\theta^d):
(\partial^\theta)^\alpha x\in
L_p(\mathbb{T}_\theta^d)
\text{ for }|\alpha|\leq k
\right\}.
\end{equation}
We equip it with the norm
\begin{equation*}
\|x\|_{W_p^k}
=\sum_{|\alpha|\leq k}
\|(\partial^\theta)^\alpha x\|_p.
\end{equation*}
For comparison, the potential-Sobolev space is
$H_p^s(\mathbb{T}_\theta^d)
=\{x\in\mathcal S':J_\theta^s x\in L_p\}$,
with norm $\|J_\theta^s x\|_p$.
We use $W_p^k$ in the endpoint embeddings,
including $p=1$; no identification with
$H_p^k$ at $p=1$ is assumed. These
conventions agree with
\cite[Section~2]{XiongXuYin} and the
differential setup of \cite[Section~2.6]{STTZ}.

Fix $\chi\in C_c^\infty(\mathbb{R}^d)$ equal
to $1$ on $|\xi|\leq1$ and supported in
$|\xi|<3/2$, and set
$\varphi(\xi)=\chi(\xi)-\chi(2\xi)$.
Then $\varphi$ is supported in an annulus and
\begin{equation*}
\chi(\xi)+\sum_{j\geq1}
\varphi(2^{-j}\xi)=1
\ \ (\xi\in\mathbb{R}^d).
\end{equation*}
For a distribution $x$, set
$\Delta_0x=T_{\chi|_{\mathbb{Z}^d}}^\theta x$
and $\Delta_jx=T_{\varphi(2^{-j}\cdot)|_{\mathbb{Z}^d}}^\theta x$
for $j\geq1$, using the multiplier convention
\eqref{eq:prelim-multiplier}. Equivalently,
$\Delta_0=\chi(D)$ and
$\Delta_j=\varphi(2^{-j}D)$, where
$D=(-i\partial_1^\theta,\ldots,-i\partial_d^\theta)$.
For $s\in\mathbb{R}$ and $1\leq p,q\leq\infty$,
define the inhomogeneous Besov space
$B^s_{p,q}(\mathbb{T}_\theta^d)$ as the set of
$x\in\mathcal S'(\mathbb{T}_\theta^d)$ for which
\begin{equation}\label{eq:prelim-Besov}
\|x\|_{B^s_{p,q}(\mathbb{T}_\theta^d)}
=\left\|
\bigl(2^{js}\|\Delta_jx\|_{L_p(\mathbb{T}_\theta^d)}
\bigr)_{j\geq0}
\right\|_{\ell_q}<\infty.
\end{equation}
Here the $\ell_\infty$ norm is the supremum.

For the target $B_{\infty,\infty}^s$, we fix
the following equivalent mean-plus-annular-block
norm. For every $j\geq0$, let
\begin{equation*}
\widetilde\varphi_j(t)
=\sum_{n\in\mathbb{Z}^d}
\varphi(2^{-j}n)e^{in\cdot t},
\ \  t\in\mathbb{T}^d.
\end{equation*}
This is a finite trigonometric polynomial with
$\widehat{\widetilde\varphi_j}(n)
=\varphi(2^{-j}n)$. Its semi-convolution
action extends coefficientwise to distributions:
\begin{equation*}
\widetilde\varphi_j\ast x
:=T_{\varphi(2^{-j}\cdot)|_{\mathbb{Z}^d}}^\theta x
=\sum_{n\in\mathbb{Z}^d}
\varphi(2^{-j}n)\widehat x(n)e_n^\theta.
\end{equation*}
For $s\in\mathbb{R}$, the norm used in the
endpoint embeddings is
\begin{equation}\label{eq:Besov-infinity-norm}
\|x\|_{B_{\infty,\infty}^s(\mathbb{T}_\theta^d)}
:=\max\left\{
|\widehat x(0)|,\
\sup_{j\geq0}2^{js}
\|\widetilde\varphi_j\ast x\|_{L_\infty(\mathbb{T}_\theta^d)}
\right\}.
\end{equation}
Indeed, $\chi(2n)=0$ for every nonzero
$n\in\mathbb{Z}^d$, whence
$\Delta_0x=\widehat x(0)\mathbf1+
\widetilde\varphi_0\ast x$ and
$\Delta_jx=\widetilde\varphi_j\ast x$ for $j\geq1$.
Thus \eqref{eq:Besov-infinity-norm} and
\eqref{eq:prelim-Besov} with $p=q=\infty$
define the same space and equivalent norms.
Moreover,
$\sum_{j\geq0}\varphi(2^{-j}n)=1$
for $n\in\mathbb{Z}^d\setminus\{0\}$.
Different admissible dyadic partitions yield
equivalent norms
\cite[Definition~3.1 and Remark~3.2]{XiongXuYin}.

\section{Convolution and sharp endpoint bounds}
\label{sec:convolution}

We use the Fourier realization of noncommutative tori described in \cite[Section~1]{ChenXuYin} and \cite[Section~1.2]{XiongXuYin}.

We write
\begin{equation*}
S(\mathbb{T}_{\theta}^{d}) = S\bigl(L_{\infty}(\mathbb{T}_{\theta}^{d}),\tau_{\theta}\bigr)
\end{equation*}
for the algebra of $\tau_{\theta}$-measurable operators affiliated with $L_{\infty}(\mathbb{T}_{\theta}^{d})$.

\subsection{Construction and algebraic properties}

For a deformation parameter $\theta$, let $\{e_{n}^{\theta}:n\in\mathbb{Z}^{d}\}$ be the Weyl Fourier system.  Thus
\begin{equation*}
e_{n}^{\theta}e_{m}^{\theta} = \sigma_{\theta}(n,m)e_{n+m}^{\theta},
\ \  (e_{n}^{\theta})^{*}=e_{-n}^{\theta},
\end{equation*}
where
\begin{equation*}
\sigma_{\theta}(n,m) = \exp\bigl(\pi i\langle n,\theta m\rangle\bigr).
\end{equation*}
The canonical trace is determined by
\begin{equation}\label{eq:canonical-trace-basis}
\tau_{\theta}(e_{n}^{\theta}) = 
\begin{cases}
1, & n=0,\\
0, & n\neq0.
\end{cases}
\end{equation}

Equivalently, $\tau_{\theta}(e_{n}^{\theta})=\delta_{n,0}$, where $\delta_{n,0}$ is the Kronecker delta.  Since $\theta$ is skew-symmetric, $\langle n,\theta n\rangle=0$, and therefore
\begin{equation*}
\tau_{\theta}\bigl((e_{n}^{\theta})^{*}e_{m}^{\theta}\bigr) = \delta_{n,m}.
\end{equation*}
For $x\in L_{1}(\mathbb{T}_{\theta}^{d})$, we define
\begin{equation*}
\widehat{x}(n) = \tau_{\theta}\bigl(x(e_{n}^{\theta})^{*}\bigr) = \tau_{\theta}(xe_{-n}^{\theta}).
\end{equation*}
Our derivations satisfy $\partial_{j} e_n^\theta=in_{j} e_n^\theta$. The basic cocycle identity used below is
\begin{equation}\label{eq:cocycle-addition}
\sigma_{\theta_{1}+\theta_{2}}(n,m) = \sigma_{\theta_{1}}(n,m)\sigma_{\theta_{2}}(n,m).
\end{equation}

\begin{proposition}\label{prop:comultiplication}
For all deformation parameters $\theta_{1}$ and $\theta_{2}$, there exists a unique normal unital injective $*$-homomorphism
\begin{equation*}
\mathcal{T}_{\theta_{1},\theta_{2}}: L_{\infty}(\mathbb{T}_{\theta_{1}+\theta_{2}}^{d}) \to 
L_{\infty}(\mathbb{T}_{\theta_{1}}^{d}) \overline\otimes L_{\infty}(\mathbb{T}_{\theta_{2}}^{d})
\end{equation*}
satisfying
\begin{equation}\label{eq:comultiplication-basis}
\mathcal{T}_{\theta_{1},\theta_{2}} (e_{n}^{\theta_{1}+\theta_{2}}) = e_{n}^{\theta_{1}}\otimes e_{n}^{\theta_{2}},
\ \  n\in\mathbb{Z}^{d}.
\end{equation}

Moreover, it preserves the canonical traces:
\begin{equation}\label{eq:comultiplication-trace}
(\tau_{\theta_{1}}\otimes\tau_{\theta_{2}}) \circ\mathcal{T}_{\theta_{1},\theta_{2}} = \tau_{\theta_{1}+\theta_{2}}.
\end{equation}
\end{proposition}

\begin{proof}
Put
\begin{equation*}
v_{n}=e_{n}^{\theta_{1}}\otimes e_{n}^{\theta_{2}}, \ \  n\in\mathbb{Z}^{d}.
\end{equation*}
By \eqref{eq:cocycle-addition},
\begin{align*}
v_{n}v_m &= \sigma_{\theta_{1}}(n,m)\sigma_{\theta_{2}}(n,m) e_{n+m}^{\theta_{1}}\otimes e_{n+m}^{\theta_{2}} \\
&= \sigma_{\theta_{1}+\theta_{2}}(n,m)v_{n+m},
\end{align*}
and $v_{n}^{*}=v_{-n}$. Thus $(v_{n})_{n\in\mathbb{Z}^{d}}$ is a projective unitary representation with cocycle $\sigma_{\theta_{1}+\theta_{2}}$. The universal property of the twisted group $C^{*}$-algebra gives a unital $*$-homomorphism
\begin{equation}\label{eq:comultiplication-Cstar}
\mathcal T_{\theta_{1},\theta_{2}}: C(\mathbb{T}_{\theta_{1}+\theta_{2}}^{d})
\to  C(\mathbb{T}_{\theta_{1}}^{d})\otimes_{\min} C(\mathbb{T}_{\theta_{2}}^{d})
\end{equation}
satisfying \eqref{eq:comultiplication-basis}.  Here the minimal tensor product is the spatial $C^{*}$-tensor product represented inside the corresponding von Neumann tensor product.  Compare \cite[Definition~3.13 and Remark~3.8]{LeeXiong} and \cite[proof of Proposition~2.1]{Ricard}.

For every $n\in\mathbb{Z}^{d}$,
\begin{align*}
(\tau_{\theta_{1}}\otimes\tau_{\theta_{2}}) \bigl(\mathcal T_{\theta_{1},\theta_{2}} (e_{n}^{\theta_{1}+\theta_{2}})\bigr)
& = \tau_{\theta_{1}}(e_{n}^{\theta_{1}}) \tau_{\theta_{2}}(e_{n}^{\theta_{2}}) \\
& = \delta_{n,0}
= \tau_{\theta_{1}+\theta_{2}}(e_{n}^{\theta_{1}+\theta_{2}}).
\end{align*}
Linearity, norm continuity and density of the Fourier polynomials extend this trace identity to $C(\mathbb{T}_{\theta_{1}+\theta_{2}}^{d})$. In particular, the map in \eqref{eq:comultiplication-Cstar} is injective, because the canonical trace is faithful.

Apply the normal-extension argument recorded in Lemma~\ref{lem:trace-preserving-normal-extension} to this trace-preserving $*$-isomorphism onto its range and then include the von Neumann closure of that range into
\begin{equation*}
L_{\infty}(\mathbb{T}_{\theta_{1}}^{d})
\,\overline\otimes\,
L_{\infty}(\mathbb{T}_{\theta_{2}}^{d}).
\end{equation*}
This gives a normal unital injective $*$-homomorphism satisfying \eqref{eq:comultiplication-basis} and \eqref{eq:comultiplication-trace}.  Uniqueness follows because Fourier polynomials are weak-$*$ dense and normal maps are weak-$*$ continuous.
\end{proof}

\begin{remark}\label{rem:twisted-Fourier-comultiplication}
Under the canonical trace-preserving identification
\begin{equation*}
L_{\infty}(\mathbb{T}_{\theta}^{d}) \cong VN(\mathbb{Z}^{d},\sigma_{\theta}),
\end{equation*}
one has
\begin{equation*}
e_{n}^{\theta} \leftrightarrow \lambda_{\sigma_{\theta}}(n), \ \  n\in\mathbb{Z}^{d}.
\end{equation*}
See \cite[Section~7.1]{ArhancetSpectral} and \cite[Definition~3.15 and Example~3.16 (1)]{LeeXiong}.

This realization is also related to the cocycle-twisted eigenunitary and crossed-product framework considered in \cite[Section~4]{DoddsFerlegerDePagterSukochev}.  In particular, the projective unitaries appearing there satisfy
\begin{equation*}
U_\gamma U_\mu=\beta(\gamma,\mu)U_{\gamma+\mu};
\end{equation*}
see \cite[(4.1)]{DoddsFerlegerDePagterSukochev}.  For $H=\mathbb{Z}^{d}$ and $\beta=\sigma_{\theta}$, this is precisely the multiplication law of the Fourier generators $e_{n}^{\theta}$.

Since
\begin{equation*}
\sigma_{\theta_{1}+\theta_{2}} =\sigma_{\theta_{1}}\sigma_{\theta_{2}},
\end{equation*}
Proposition~\ref{prop:comultiplication} becomes
\begin{equation*}
VN(\mathbb{Z}^{d},\sigma_{\theta_{1}}\sigma_{\theta_{2}})
\rightarrow
VN(\mathbb{Z}^{d},\sigma_{\theta_{1}})
\,\overline\otimes\,
VN(\mathbb{Z}^{d},\sigma_{\theta_{2}}),
\end{equation*}
determined by
\begin{equation*}
\lambda_{\sigma_{\theta_{1}}\sigma_{\theta_{2}}}(n) \mapsto \lambda_{\sigma_{\theta_{1}}}(n) \otimes \lambda_{\sigma_{\theta_{2}}}(n).
\end{equation*}
The corresponding parameter-addition embedding for quantum tori is used explicitly in \cite[proof of Proposition~2.1]{Ricard}.

More generally, for a discrete group $G$ and scalar $2$-cocycles $\sigma_{1},\sigma_{2}$, the tensor product
\begin{equation*}
s\mapsto \lambda_{\sigma_{1}}(s)\otimes\lambda_{\sigma_{2}}(s)
\end{equation*}
is a projective unitary representation with cocycle $\sigma_{1}\sigma_{2}$; see \cite[Remark~3.8]{LeeXiong}. 
\end{remark}

\begin{proposition}\label{prop:coassociativity}
The family $\{\mathcal{T}_{\theta_{1},\theta_{2}}\}$ is coassociative.  More precisely,
\begin{align*}
(\mathcal{T}_{\theta_{1},\theta_{2}}\otimes\operatorname{id}) \circ \mathcal{T}_{\theta_{1}+\theta_{2},\theta_{3}} 
= (\operatorname{id}\otimes\mathcal{T}_{\theta_{2},\theta_{3}}) \circ \mathcal{T}_{\theta_{1},\theta_{2}+\theta_{3}}.
\end{align*}
\end{proposition}

\begin{proof}
Both sides are normal maps and send $e_{n}^{\theta_{1}+\theta_{2}+\theta_{3}}$ to
\begin{equation*}
e_{n}^{\theta_{1}} \otimes e_{n}^{\theta_{2}} \otimes e_{n}^{\theta_{3}}.
\end{equation*}
They therefore agree on Fourier polynomials and, by normality, on the whole von Neumann algebra.
\end{proof}

\subsubsection{Definition and elementary properties}

For $1\leq p\leq\infty$, we write
\begin{equation*}
L_{p}^{+}(\mathbb{T}_{\theta}^{d}) = \{x\in L_{p}(\mathbb{T}_{\theta}^{d}):x\geq0\}.
\end{equation*}
Here, $x\geq0$ means that $x$ is self-adjoint and its spectral measure is supported in $[0,\infty)$.  In particular, if $x\in L_{1}^{+}(\mathbb{T}_{\theta}^{d})$, then $\|x\|_{1}=\tau_{\theta}(x).$ For a self-adjoint measurable operator $x$, we use the standard positive and negative parts
\begin{equation*}
x=x_{+}-x_{-}, \ \  x_{\pm}=\frac{|x|\pm x}{2}.
\end{equation*}
Then $x_{+},x_{-}\geq0$, $x_{+}x_{-}=0$, and
\begin{equation*}
|x|=x_{+}+x_{-}.
\end{equation*}
The predual of $L_{\infty}(\mathbb{T}_{\theta}^{d})$ is identified with $L_{1}(\mathbb{T}_{\theta}^{d})$ through
\begin{equation*}
\langle x,a\rangle = \tau_{\theta}(xa).
\end{equation*}
Since $\mathcal{T}_{\theta_{1},\theta_{2}}$ is normal, it has a bounded preadjoint; see \cite[Chapter~III, Section~2]{Takesaki}.  For $x\in L_{1}(\mathbb{T}_{\theta_{1}}^{d})$ and $y\in L_{1}(\mathbb{T}_{\theta_{2}}^{d})$, the elementary tensor $x\otimes y$ belongs to the $L_{1}$-space of the tensor-product algebra and
\begin{equation}\label{eq:tensor-L1-norm}
\|x\otimes y\|_{1} = \|x\|_{1}\|y\|_{1}.
\end{equation}
This follows from the polar decompositions of $x$ and $y$ and the product trace.
\begin{definition}\label{def:convolution}
For $x\in L_{1}(\mathbb{T}_{\theta_{1}}^{d})$ and $y\in L_{1}(\mathbb{T}_{\theta_{2}}^{d})$, define their \emph{convolution} by
\begin{equation}\label{eq:convolution-definition}
x\star_{\theta_{1},\theta_{2}}y = (\mathcal{T}_{\theta_{1},\theta_{2}})_{*}(x\otimes y) \in L_{1}(\mathbb{T}_{\theta_{1}+\theta_{2}}^{d}).
\end{equation}

Equivalently, it is the unique element satisfying
\begin{equation}\label{eq:convolution-duality}
\tau_{\theta_{1}+\theta_{2}} \bigl((x\star_{\theta_{1},\theta_{2}}y)a\bigr)  
= (\tau_{\theta_{1}}\otimes\tau_{\theta_{2}}) \bigl((x\otimes y)\mathcal{T}_{\theta_{1},\theta_{2}}(a)\bigr)   
\end{equation}
for every $a\in L_{\infty}(\mathbb{T}_{\theta_{1}+\theta_{2}}^{d})$.
\end{definition}

The right-hand side of \eqref{eq:convolution-duality} is a normal functional of $a$, and
\begin{align*}
\left| (\tau_{\theta_{1}}\otimes\tau_{\theta_{2}}) \bigl((x\otimes y)\mathcal{T}_{\theta_{1},\theta_{2}}(a)\bigr) \right|
\leq  \|x\|_{1}\|y\|_{1}\|a\|_{\infty}.
\end{align*}
Thus the definition is well posed.  It also makes bilinearity immediate: the map $(x,y)\mapsto x\otimes y$ is bilinear and the preadjoint $(\mathcal{T}_{\theta_{1},\theta_{2}})_{*}$ is linear.

Since $\mathcal{T}_{\theta_{1},\theta_{2}}$ is unital,
\begin{equation*}
\|(\mathcal{T}_{\theta_{1},\theta_{2}})_{*}\| = \|\mathcal{T}_{\theta_{1},\theta_{2}}\| =1.
\end{equation*}
The preadjoint is positive.  Indeed, if $h\geq0$ belongs to the $L_{1}$-space of the tensor-product algebra and $a\geq0$, then
\begin{align*}
\left\langle (\mathcal{T}_{\theta_{1},\theta_{2}})_{*}(h),a \right\rangle
&= (\tau_{\theta_{1}}\otimes\tau_{\theta_{2}}) \bigl(h\mathcal{T}_{\theta_{1},\theta_{2}}(a)\bigr) \\
&= (\tau_{\theta_{1}}\otimes\tau_{\theta_{2}}) \bigl(h^{\frac{1}{2}}\mathcal{T}_{\theta_{1},\theta_{2}}(a)h^{\frac{1}{2}}\bigr) \geq0.
\end{align*}

\begin{proposition}\label{prop:basic-properties}
Let $x\in L_{1}(\mathbb{T}_{\theta_{1}}^{d})$ and $y\in L_{1}(\mathbb{T}_{\theta_{2}}^{d})$. Then:

\begin{enumerate}
\item $\widehat{x\star_{\theta_{1},\theta_{2}}y}(n) =\widehat{x}(n)\widehat{y}(n)$ for every $n\in\mathbb{Z}^{d}$;

\item $\tau_{\theta_{1}+\theta_{2}} (x\star_{\theta_{1},\theta_{2}}y) =\tau_{\theta_{1}}(x)\tau_{\theta_{2}}(y)$;

\item if $x,y\geq0$, then $x\star_{\theta_{1},\theta_{2}}y\geq0$;

\item 
\begin{equation}\label{eq:L1}
\|x\star_{\theta_{1},\theta_{2}}y\|_{1} \leq \|x\|_{1}\|y\|_{1}.
\end{equation}
\end{enumerate}
If $x,y\geq0$, equality holds in \eqref{eq:L1}.
\end{proposition}

\begin{proof}
Put $z=x\star_{\theta_{1},\theta_{2}}y$.  Taking $a=e_{-n}^{\theta_{1}+\theta_{2}}$ in \eqref{eq:convolution-duality} gives
\begin{align*}
\widehat{z}(n) &= (\tau_{\theta_{1}}\otimes\tau_{\theta_{2}}) \bigl((x\otimes y)
(e_{-n}^{\theta_{1}}\otimes e_{-n}^{\theta_{2}})\bigr)
= \widehat{x}(n)\widehat{y}(n).
\end{align*}
Taking $a=\1$ gives the trace formula. Suppose that $x,y\geq0$ and put $h=x\otimes y$.  If $a\geq0$, then $\mathcal{T}_{\theta_{1},\theta_{2}}(a)\geq0$, and traciality gives
\begin{align*}
\tau_{\theta_{1}+\theta_{2}}(za) &= (\tau_{\theta_{1}}\otimes\tau_{\theta_{2}})
\bigl(h\mathcal{T}_{\theta_{1},\theta_{2}}(a)\bigr) \\
&= (\tau_{\theta_{1}}\otimes\tau_{\theta_{2}}) \bigl(h^{\frac{1}{2}}\mathcal{T}_{\theta_{1},\theta_{2}}(a)h^{\frac{1}{2}}\bigr) \geq0.
\end{align*}
Hence $z\geq0$.

Finally, $L_{1}$--$L_{\infty}$ duality, \eqref{eq:tensor-L1-norm}, and contractivity of the unital $*$-homomorphism $\mathcal{T}_{\theta_{1},\theta_{2}}$ give
\begin{align*}
\|z\|_{1} & = \sup_{\|a\|_{\infty}\leq1} \left|\tau_{\theta_{1}+\theta_{2}}(za)\right| \leq \|x\otimes y\|_{1} \sup_{\|a\|_{\infty}\leq1} \|\mathcal{T}_{\theta_{1},\theta_{2}}(a)\|_{\infty}
\leq \|x\|_{1}\|y\|_{1}.
\end{align*}
When $x,y\geq0$, positivity and the trace formula imply
\begin{equation*}
\|z\|_{1} = \tau_{\theta_{1}+\theta_{2}}(z) = \|x\|_{1}\|y\|_{1}.
\end{equation*}
\end{proof}

For Fourier polynomials
\begin{equation*}
x=\sum_{n}a_{n}e_{n}^{\theta_{1}}, \ \  y=\sum_{n}b_{n}e_{n}^{\theta_{2}},
\end{equation*}
we consequently have
\begin{equation}\label{eq:polynomial-formula}
x\star_{\theta_{1},\theta_{2}}y = \sum_{n}a_{n}b_{n}e_{n}^{\theta_{1}+\theta_{2}}.
\end{equation}

Thus convolution is coefficientwise multiplication of Fourier coefficients; it is not the twisted product inside a fixed noncommutative torus.

\begin{proposition}\label{prop:associative-commutative} 
The family of convolutions is associative and commutative under exchange of the inputs and their deformation parameters. Namely,
\begin{align*}
(x\star_{\theta_{1},\theta_{2}}y) \star_{\theta_{1}+\theta_{2},\theta_{3}}z  =
x\star_{\theta_{1},\theta_{2}+\theta_{3}} (y\star_{\theta_{2},\theta_{3}}z)
\end{align*}
and
\begin{equation*}
x\star_{\theta_{1},\theta_{2}}y = y\star_{\theta_{2},\theta_{1}}x.
\end{equation*}
\end{proposition}

\begin{proof}
Both sides of the associativity identity have $n$-th Fourier coefficient
\begin{equation*}
\widehat{x}(n)\widehat{y}(n)\widehat{z}(n).
\end{equation*}
Both sides of the commutativity identity have $n$-th Fourier coefficient
\begin{equation*}
\widehat{x}(n)\widehat{y}(n).
\end{equation*}
Fourier uniqueness in $L_{1}(\mathbb{T}_{\theta}^{d})$ gives both assertions.  Indeed, if all Fourier coefficients of an $L_{1}$-operator vanish, then its trace pairing with every Fourier polynomial vanishes; weak-$*$ density of Fourier polynomials in $L_{\infty}(\mathbb{T}_{\theta}^{d})$ then implies that the operator is zero. See also \cite[Section~1.2]{XiongXuYin}.
\end{proof}

\subsubsection{The zero-parameter case: semi-convolution}

Let $L_{\infty}(\mathbb{T}_{0}^{d})=L_{\infty}(\mathbb{T}^{d})$, and let $m$ be normalized Haar measure.  The canonical action is determined by
\begin{equation*}
\alpha_{s}(e_{n}^{\theta}) = e^{in\cdot s}e_{n}^{\theta}.
\end{equation*}
Each $\alpha_s$ is a trace-preserving $*$-automorphism of $L_{\infty}(\mathbb{T}_{\theta}^{d})$. For $1\leq p<\infty$, its extension to $L_{p}$ belongs to a strongly continuous group of isometries; see \cite[Lemma~1.1]{XiongXuYin}. For $f\in L_{1}(\mathbb{T}^{d})$ and $x\in L_{1}(\mathbb{T}_{\theta}^{d})$, define the semi-convolution
\begin{equation}\label{eq:semi-convolution-definition}
f\ast x = \int\limits_{\mathbb{T}^{d}}f(s)\alpha_{-s}(x)\,dm(s),
\end{equation}
where the integral is a Bochner integral in $L_{1}(\mathbb{T}_{\theta}^{d})$.  With
\begin{equation*}
\widehat{f}(n) = \int\limits_{\mathbb{T}^{d}}f(s)e^{-in\cdot s}\,dm(s),
\end{equation*}
we have
\begin{equation*}
\widehat{f\ast x}(n) = \widehat{f}(n)\widehat{x}(n).
\end{equation*}

Proposition~\ref{prop:basic-properties} gives the same Fourier coefficients for $f\star_{0,\theta}x$.  Fourier uniqueness therefore yields the exact identity
\begin{equation}\label{eq:semi-is-convolution}
f\star_{0,\theta}x = f\ast x.
\end{equation}
If semi-convolution is defined using $\alpha_{s}$ instead of $\alpha_{-s}$, then \eqref{eq:semi-is-convolution} holds with $f$ replaced by $s\mapsto f(-s)$; this reflection does not change its decreasing rearrangement.

\subsection{Sharp endpoint estimates}\label{subsec:endpoints}

We construct the parameter-reversing map used in the sharp mixed endpoint estimate.  To make the construction transparent, we separate the algebraic, $C^{*}$-algebraic, von Neumann algebraic, $L_{p}$ and Fourier steps.

\begin{lem}\label{lem:parameter-reversal}
For every deformation parameter $\theta$, there exists a unique trace-preserving normal linear $*$-anti-isomorphism
\begin{equation*}
\kappa_{\theta}: L_{\infty}(\mathbb{T}_{\theta}^{d}) \to  L_{\infty}(\mathbb{T}_{-\theta}^{d})
\end{equation*}
satisfying
\begin{equation}\label{eq:kappa-basis}
\kappa_{\theta}(e_{n}^{\theta})=e_{-n}^{-\theta}, \ \  n\in\mathbb{Z}^{d}.
\end{equation}
It extends to a surjective linear isometry between the corresponding $L_{p}$-spaces for every $1\leq p<\infty$, and to a map
\begin{equation*}
\kappa_{\theta}: S(\mathbb{T}_{\theta}^{d}) \to  S(\mathbb{T}_{-\theta}^{d})
\end{equation*}
on the algebras of measurable operators.  For every $x\in S(\mathbb{T}_{\theta}^{d})$,

\begin{equation}\label{eq:kappa-mu}
\mu(s,\kappa_{\theta}(x))=\mu(s,x), \ \  0<s<1,
\end{equation}
and, for every $x\in L_{1}(\mathbb{T}_{\theta}^{d})$,
\begin{equation}\label{eq:kappa-Fourier}
\widehat{\kappa_{\theta}(x)}(n)=\widehat{x}(-n), \ \  n\in\mathbb{Z}^{d}.
\end{equation}
\end{lem}

The construction and proof of Lemma~\ref{lem:parameter-reversal} are given in Appendix~\ref{app:parameter-reversal}.

\begin{proposition}\label{prop:L2-endpoint}
If $x\in L_{2}(\mathbb{T}_{\theta_{1}}^{d})$ and $y\in L_{2}(\mathbb{T}_{\theta_{2}}^{d})$, then
\begin{equation*}
\|x\star_{\theta_{1},\theta_{2}}y\|_{\infty}
\leq \|x\|_{2}\|y\|_{2}.
\end{equation*}
\end{proposition}

\begin{proof}
Cauchy--Schwarz and Plancherel give
\begin{equation*}
\sum_{n\in\mathbb{Z}^{d}} |\widehat{x}(n)\widehat{y}(n)| \leq \|x\|_{2}\|y\|_{2}.
\end{equation*}
Hence
\begin{equation*}
\sum_{n\in\mathbb{Z}^{d}} \widehat{x}(n)\widehat{y}(n)e_{n}^{\theta_{1}+\theta_{2}}
\end{equation*}
converges absolutely in operator norm, with norm at most $\|x\|_{2}\|y\|_{2}$.  Since the traces are finite, $x$ and $y$ also belong to $L_{1}$.  Proposition~\ref{prop:basic-properties} and Fourier uniqueness identify the sum with $x\star_{\theta_{1},\theta_{2}}y$.
\end{proof}

\begin{lem}[Preadjoint representation of convolution]
\label{lem:Sx-preadjoint}
Let $\theta_{1}$ and $\theta_{2}$ be deformation parameters and put $\theta_{3}=\theta_{1}+\theta_{2}$. For $x\in L_{1}(\mathbb{T}_{\theta_{1}}^{d})$, define
\begin{equation}\label{eq:def-Sx}
S_{x}:L_{1}(\mathbb{T}_{\theta_{3}}^{d}) \to  L_{1}(\mathbb{T}_{\theta_{2}}^{d}), \ \  
S_{x}(a)=\kappa_{\theta_{1}}(x) \star_{-\theta_{1},\theta_{3}}a.
\end{equation}
Then
\begin{equation}\label{eq:Sx-norm}
\|S_{x}\|_{L_{1}(\mathbb{T}_{\theta_3}^{d}) \to L_{1}(\mathbb{T}_{\theta_{2}}^{d})}\leq\|x\|_{1}.
\end{equation}
The adjoint map
\begin{equation*}
S_{x}^{*}:L_\infty(\mathbb{T}_{\theta_{2}}^{d}) \to  L_\infty(\mathbb{T}_{\theta_3}^{d})
\end{equation*}
is normal, $S_{x}$ is its preadjoint, and
\begin{equation}\label{eq:Sx-adjoint-convolution}
S_{x}^{*}(y)=x\star_{\theta_{1},\theta_{2}}y,
\ \  y\in L_\infty(\mathbb{T}_{\theta_{2}}^{d}).
\end{equation}
Consequently,
\begin{equation}\label{eq:Sx-adjoint-norm}
\|x\star_{\theta_{1},\theta_{2}}y\|_\infty
\leq\|x\|_{1}\|y\|_\infty.
\end{equation}
\end{lem}

\begin{proof}
The deformation parameter of the range in \eqref{eq:def-Sx} is $-\theta_{1}+\theta_3=\theta_{2}$. Proposition~\ref{prop:basic-properties} and the $L_{1}$-isometry of $\kappa_{\theta_{1}}$ give
\begin{equation*}
\|S_{x}(a)\|_{1} \leq \|\kappa_{\theta_{1}}(x)\|_{1}\|a\|_{1}
=\|x\|_{1}\|a\|_{1},
\end{equation*}
which proves \eqref{eq:Sx-norm}. Its Banach adjoint satisfies
\begin{equation}\label{eq:Sx-adjoint}
\tau_{\theta_3}\bigl((S_{x}^{*}y)a\bigr) = \tau_{\theta_{2}}\bigl(yS_{x}(a)\bigr)
\end{equation}
for $a\in L_{1}(\mathbb{T}_{\theta_3}^{d})$ and $y\in L_\infty(\mathbb{T}_{\theta_{2}}^{d})$. Since $S_{x}$ acts between the canonical preduals, $S_{x}^{*}$ is normal and $S_{x}$ is its preadjoint.

For $n\in\mathbb{Z}^{d}$, the Fourier formula for convolution and \eqref{eq:kappa-Fourier} give
\begin{equation*}
S_{x}(e_{-n}^{\theta_3})
=\widehat{\kappa_{\theta_{1}}(x)}(-n)e_{-n}^{\theta_{2}}
=\widehat x(n)e_{-n}^{\theta_{2}}.
\end{equation*}
Therefore \eqref{eq:Sx-adjoint} yields
\begin{equation*}
\widehat{S_{x}^{*}y}(n) = \widehat x(n)\widehat y(n).
\end{equation*}
Both $S_{x}^{*}y$ and $x\star_{\theta_{1},\theta_{2}}y$ belong to $L_{1}(\mathbb{T}_{\theta_3}^{d})$ and have the same Fourier coefficients. Fourier uniqueness proves \eqref{eq:Sx-adjoint-convolution}. Finally, $\|S_{x}^{*}\|=\|S_{x}\|\leq\|x\|_{1}$, which gives \eqref{eq:Sx-adjoint-norm}.
\end{proof}

\begin{proposition}\label{prop:L1-Linfty}
For all deformation parameters $\theta_{1}$ and $\theta_{2}$,
\begin{equation}\label{eq:L1-Linfty-sharp}
\|x\star_{\theta_{1},\theta_{2}}y\|_{\infty} \leq \|x\|_{1}\|y\|_{\infty}
\end{equation}
for every $x\in L_{1}(\mathbb{T}_{\theta_{1}}^{d})$ and $y\in L_{\infty}(\mathbb{T}_{\theta_{2}}^{d})$.  Similarly, 
\begin{equation}\label{eq:Linfty-L1-sharp}
\|x\star_{\theta_{1},\theta_{2}}y\|_{\infty} \leq \|x\|_{\infty}\|y\|_{1}
\end{equation}
for every $x\in L_{\infty}(\mathbb{T}_{\theta_{1}}^{d})$ and $y\in L_{1}(\mathbb{T}_{\theta_{2}}^{d})$.  Both constants are optimal.
\end{proposition}

\begin{proof}
The first estimate is \eqref{eq:Sx-adjoint-norm}.  Interchanging the two inputs and the two deformation parameters gives the second estimate. Finally, taking $x=y=\1$ shows that neither constant can be smaller than one.
\end{proof}

\begin{remark}[Opposite deformations and scalar-valued convolution]
\label{rem:diamond-convolution}
The reflection $\rho_\theta(e_n^\theta)=e_{-n}^\theta$ is a trace-preserving automorphism. Thus $\iota_\theta=\kappa_\theta\circ\rho_\theta$ is a trace-preserving anti-isomorphism from the parameter-$\theta$ torus to the parameter-$-\theta$ torus, with $\widehat{\iota_\theta(y)}(n)=\widehat y(n)$ and $\mu(\iota_\theta(y))=\mu(y)$. Consequently,
\begin{equation*}
x\diamond_\theta y := x\star_{\theta,-\theta}\iota_\theta(y)
\in L_{1}(\mathbb{T}^{d}), \ \  x,y\in L_{1}(\mathbb{T}_{\theta}^{d}),
\end{equation*}
has Fourier coefficients $\widehat x(n)\widehat y(n)$. It is commutative, positive for positive inputs, and satisfies $\|x\diamond_\theta y\|_{1}\leq\|x\|_{1}\|y\|_{1}$. Theorem~\ref{thm:Young} supplies its Young bounds. Related scalar-valued operator convolutions occur in quantum harmonic analysis \cite{Werner,LuefSkrettingland}.
\end{remark}

\section{Sharp Young inequalities}\label{sec:young}

The following result is the noncommutative counterpart of the classical Young convolution inequality; see \cite{Young}.

\begin{thm}[Sharp Young inequality]\label{thm:Young}
Let $1\leq p,q,r\leq\infty$ satisfy
\begin{equation*}
1+\frac{1}{r} = \frac{1}{p}+\frac{1}{q}.
\end{equation*}
Then, for every $x\in L_{p}(\mathbb{T}_{\theta_{1}}^{d})$ and $y\in L_{q}(\mathbb{T}_{\theta_{2}}^{d})$, one has
\begin{equation*}
x\star_{\theta_{1},\theta_{2}}y \in L_{r}(\mathbb{T}_{\theta_{1}+\theta_{2}}^{d})
\end{equation*}
and
\begin{equation}\label{eq:Young}
\left\|
x\star_{\theta_{1},\theta_{2}}y
\right\|_{L_{r}(\mathbb{T}_{\theta_{1}+\theta_{2}}^{d})}
\leq
\|x\|_{L_{p}(\mathbb{T}_{\theta_{1}}^{d})}
\|y\|_{L_{q}(\mathbb{T}_{\theta_{2}}^{d})}.
\end{equation}
The constant $1$ in \eqref{eq:Young} is optimal.
\end{thm}

\begin{proof}
Put $\theta_{3}=\theta_{1}+\theta_{2}$, and let $r'$ be the conjugate exponent of $r$. Since the canonical traces are finite,  the convolutions below are well defined in $L_{1}$.

For Fourier polynomials $x$, $y$, and $z$ on the corresponding quantum tori, define
\begin{equation*}
\Lambda(x,y,z) = \tau_{\theta_{3}} \left( \bigl(x\star_{\theta_{1},\theta_{2}}y\bigr)z \right).
\end{equation*}
By \eqref{eq:L1}, \eqref{eq:L1-Linfty-sharp}, and \eqref{eq:Linfty-L1-sharp},
\begin{align*}
|\Lambda(x,y,z)| & \leq \|x\|_{1}\|y\|_{1}\|z\|_{\infty},\\
|\Lambda(x,y,z)| & \leq \|x\|_{1}\|y\|_{\infty}\|z\|_{1},\\
|\Lambda(x,y,z)| & \leq \|x\|_{\infty}\|y\|_{1}\|z\|_{1}.
\end{align*}
Denote the corresponding endpoint constants by
\begin{equation*}
M_{1,1,\infty},\ \  M_{1,\infty,1},\ \  M_{\infty,1,1}.
\end{equation*}
The preceding estimates show that
\begin{equation}\label{eq:Young-endpoint-constants}
M_{1,1,\infty}, M_{1,\infty,1}, M_{\infty,1,1} \leq 1.
\end{equation}

We use multilinear complex interpolation on the common Fourier-polynomial core; see \cite[Theorems~4.4.1 and~4.4.2]{BerghLofstrom}, together with the isometric complex interpolation identities for noncommutative $L_{p}$-spaces; see \cite{Kosaki}. Recall that interpolation between bounds $M_0$ and $M_{1}$ at the parameter $\vartheta$ gives the constant $M_0^{1-\vartheta}M_{1}^\vartheta.$

Set
\begin{equation*}
\alpha=\frac{1}{r},
\ \  \beta=\frac{1}{p}-\frac{1}{r},
\ \  \gamma=\frac{1}{q}-\frac{1}{r}.
\end{equation*}
The relation between the exponents implies
\begin{equation*}
\beta=1-\frac{1}{q},
\ \  \gamma=1-\frac{1}{p},
\ \  \alpha+\beta+\gamma=1.
\end{equation*}
In particular,
\begin{align*}
\left( \frac{1}{p}, \frac{1}{q}, \frac{1}{r'} \right)
&= \alpha(1,1,0) + \beta(1,0,1) + \gamma(0,1,1).
\end{align*}

Suppose first that $1<r<\infty$. Set
\begin{equation*}
\eta=\frac{\gamma}{\beta+\gamma} =\frac{\gamma}{1-\alpha}
\end{equation*}
and define $u,v\in[1,\infty]$ by
\begin{equation*}
\frac{1}{u}=1-\eta =\frac{\beta}{1-\alpha},
\ \ 
\frac{1}{v}= \eta =\frac{\gamma}{1-\alpha}.
\end{equation*}

First, interpolate between the endpoint estimates
\begin{equation*}
L_{1}\times L_{\infty}\times L_{1}\to\mathbb C
\end{equation*}
and
\begin{equation*}
L_{\infty}\times L_{1}\times L_{1}\to\mathbb C
\end{equation*}
at the parameter $\eta$. This gives
\begin{equation*}
|\Lambda(x,y,z)| \leq M_{u,v,1}\|x\|_u\|y\|_v\|z\|_{1},
\end{equation*}
where
\begin{equation*}
M_{u,v,1} \leq M_{1,\infty,1}^{1-\eta} M_{\infty,1,1}^{\eta}.
\end{equation*}

Next, interpolate this estimate with
\begin{equation*}
L_{1}\times L_{1}\times L_{\infty}\to\mathbb C
\end{equation*}
at the parameter $\alpha$. Indeed,
\begin{equation*}
\frac{1-\alpha}{u}+\alpha=\frac{1}{p}, \ \ 
\frac{1-\alpha}{v}+\alpha=\frac{1}{q}, \ \ 
1-\alpha=\frac{1}{r'}.
\end{equation*}
Therefore,
\begin{equation*}
|\Lambda(x,y,z)| \leq M_{p,q,r'}\|x\|_{p}\|y\|_{q}\|z\|_{r'},
\end{equation*}
where
\begin{align*}
M_{p,q,r'} \leq M_{u,v,1}^{1-\alpha} M_{1,1,\infty}^{\alpha} 
\leq M_{1,1,\infty}^{\alpha} M_{1,\infty,1}^{(1-\alpha)(1-\eta)} M_{\infty,1,1}^{(1-\alpha)\eta}.
\end{align*}
Since
\begin{equation*}
(1-\alpha)(1-\eta)=\beta, \ \  (1-\alpha)\eta=\gamma,
\end{equation*}
we obtain the exact constant
\begin{equation}\label{eq:Young-interpolated-constant}
M_{p,q,r'} \leq M_{1,1,\infty}^{\alpha} M_{1,\infty,1}^{\beta}
M_{\infty,1,1}^{\gamma} \leq1.
\end{equation}

If $r=\infty$, then $\alpha=0$, and the first interpolation step alone gives the required estimate. If $r=1$, the exponent relation forces $p=q=1$, so the assertion follows from \eqref{eq:L1}. The cases $\beta=0$ or $\gamma=0$ reduce to the corresponding endpoint estimates.

For $1<r<\infty$, duality applied to the bound on the polynomial core gives an element of $L_r$ with norm at most $\|x\|_{p}\|y\|_q$. When $r=\infty$, the same argument uses $L_{1}^{*}=L_\infty$; when $r=1$, the result is already \eqref{eq:L1}. Testing against Fourier polynomials identifies this element with the convolution previously defined in $L_{1}$. For polynomial inputs this proves \eqref{eq:Young}. Approximation in $L_{p}$ and $L_q$ then proves the result for all inputs with finite exponents. If $p=\infty$ or $q=\infty$, the exponent relation forces $(p,q,r)=(\infty,1,\infty)$ or $(1,\infty,\infty)$, respectively; these are the mixed endpoint estimates \eqref{eq:L1-Linfty-sharp}--\eqref{eq:Linfty-L1-sharp}.

Finally, trace preservation gives
\begin{equation*}
\mathbf{1}\star_{\theta_{1},\theta_{2}}\mathbf{1} = \mathbf{1}.
\end{equation*}
Since the canonical traces are normalized,
\begin{equation*}
\|\mathbf{1}\|_{p} = \|\mathbf{1}\|_{q} = \|\mathbf{1}\|_{r} = 1.
\end{equation*}
Hence any constant $C$ in \eqref{eq:Young} must satisfy
\begin{equation*}
1 = \left\| \mathbf{1}\star_{\theta_{1},\theta_{2}}\mathbf{1} \right\|_{r} \leq C.
\end{equation*}
Thus the optimal constant is exactly $1$.
\end{proof}

\section{Limiting estimates and submajorization}\label{sec:limiting}

We now apply the two mixed endpoint estimates to rearrangements.

\begin{corollary}\label{cor:HLP}
If $x\in L_{1}(\mathbb{T}_{\theta_{1}}^{d})$ and $y\in L_{1}(\mathbb{T}_{\theta_{2}}^{d})$, then
\begin{equation}\label{eq:HLP-first}
\mu(x\star_{\theta_{1},\theta_{2}}y) \prec\!\prec \|x\|_{1}\mu(y)
\end{equation}
and
\begin{equation}\label{eq:HLP-second}
\mu(x\star_{\theta_{1},\theta_{2}}y) \prec\!\prec \|y\|_{1}\mu(x).
\end{equation}
\end{corollary}
\begin{proof}
Fix $x$ and define
\begin{equation*}
T_{x}(y)=x\star_{\theta_{1},\theta_{2}}y.
\end{equation*}
Propositions~\ref{prop:basic-properties} and~\ref{prop:L1-Linfty} give
\begin{equation*}
\|T_{x}\|_{L_{1}\to L_{1}} \leq \|x\|_{1},
\ \  \|T_{x}\|_{L_{\infty}\to L_{\infty}} \leq \|x\|_{1}.
\end{equation*}

The second restriction is consistent with the original $L_{1}$ definition because the trace is finite and $L_{\infty}(\mathbb{T}_{\theta_{2}}^{d})\subset L_{1}(\mathbb{T}_{\theta_{2}}^{d})$.  The endpoint estimates imply
\begin{align*}
K\bigl(t,T_{x}(y);L_{1},L_{\infty}\bigr) \leq
\|x\|_{1}K\bigl(t,y;L_{1},L_{\infty}\bigr).
\end{align*}
Using \eqref{eq:K-prelim} gives \eqref{eq:HLP-first}. Commutativity under exchange of the parameters gives \eqref{eq:HLP-second}.
\end{proof}

Equivalently, for every $0<t\leq1$,
\begin{equation}\label{eq:HLP-symmetric-integral}
    \int\limits_{0}^{t} \mu(s,x\star_{\theta_{1},\theta_{2}}y)\,ds
\leq \min\left\{
\|x\|_{1}\int\limits_{0}^{t}\mu(s,y)\,ds,
\|y\|_{1}\int\limits_{0}^{t}\mu(s,x)\,ds
\right\}.
\end{equation}

\begin{corollary}[Fully symmetric convolution modules]
\label{cor:symmetric-module}
Let $E$ be a fully symmetric Banach function space on $(0,1)$. For $x\in L_{1}(\mathbb{T}_{\theta_{1}}^{d})$ and $y\in\mathcal{E}(\mathbb{T}_{\theta_{2}}^{d})$,
\begin{equation}\label{eq:symmetric-module}
\|x\star_{\theta_{1},\theta_{2}}y\|_{\mathcal{E}(\mathbb{T}_{\theta_{1}+\theta_{2}}^{d})}
\leq \|x\|_{L_{1}(\mathbb{T}_{\theta_{1}}^{d})}\|y\|_{\mathcal{E}(\mathbb{T}_{\theta_{2}}^{d})}.
\end{equation}
The analogous inequality holds with the factors interchanged.
\end{corollary}
\begin{proof}
On a probability space $\mathcal{E}\subset L_{1}$, so the convolution is defined. Apply Corollary~\ref{cor:HLP} and full symmetry of $E$.
\end{proof}

\section{Endpoint Sobolev embeddings}\label{sec:endpoint-sobolev}
\label{subsec:multiplier-reconstruction}

We begin by recording the exact compatibility between convolution and Fourier multipliers. The two endpoint embeddings below answer the questions in \cite[Remark~6.8]{XiongXuYin}. We combine Fourier reconstruction formulas with the endpoint and distribution estimates for convolution. We then apply this principle to Sobolev--Besov and critical $L_{\infty}$ embeddings.

Throughout this section, $\varphi$, $\widetilde{\varphi}_j$ and the norm \eqref{eq:Besov-infinity-norm} are those fixed in the preliminaries. In particular, 
\begin{equation*}
    \sum_{j\geq0}\varphi(2^{-j}n)=1 \textnormal{ for }n\in\mathbb{Z}^d\setminus\{0\}
\end{equation*}
and $\widetilde{\varphi}_j\ast x=\Delta_jx$.

For the remainder of this subsection, put
\begin{equation*}
D_{\ell}=-i\partial_{\ell}.
\end{equation*}
Then
\begin{equation*}
D_{\ell}e_n^{\theta}=n_{\ell}e_n^{\theta}.
\end{equation*}
Thus $D_{\ell}^{k}$ and $\partial_{\ell}^{k}$ differ only by a scalar of modulus one and have the same $L_{p}$-norm.

\begin{proposition}[convolution and Fourier multipliers]
\label{prop:multiplier-compatibility}
Let $\theta_{1}$ and $\theta_{2}$ be deformation parameters
and put $\theta_3=\theta_{1}+\theta_{2}$. We use the Fourier
multiplier $T_m^\vartheta$ on $\mathcal S'(\mathbb{T}_\vartheta^d)$
defined in \eqref{eq:prelim-multiplier} for every polynomially
bounded symbol $m:\mathbb{Z}^d\to\mathbb C$.
The convolution has a canonical extension 
\begin{equation*}
\mathcal{S}'(\mathbb{T}_{\theta_{1}}^{d})
\times
\mathcal{S}'(\mathbb{T}_{\theta_{2}}^{d})
\to
\mathcal{S}'(\mathbb{T}_{\theta_3}^{d})
\end{equation*}
determined by
\begin{equation}\label{eq:distributional-convolution}
\widehat{x\star_{\theta_{1},\theta_{2}}y}(n) = \widehat{x}(n)\widehat{y}(n).
\end{equation}
This extension agrees with Definition~\ref{def:convolution} whenever $x$ and $y$ belong to the corresponding $L_{1}$-spaces.

If $m_{1}$ and $m_{2}$ are polynomially bounded symbols, then
\begin{equation}\label{eq:two-multiplier-convolution-identity}
\bigl(T_{m_{1}}^{\theta_{1}}x\bigr) \star_{\theta_{1},\theta_{2}} \bigl(T_{m_{2}}^{\theta_{2}}y\bigr) 
= T_{m_{1}m_{2}}^{\theta_3} \bigl(x\star_{\theta_{1},\theta_{2}}y\bigr).
\end{equation}
In particular,
\begin{equation}
\label{eq:multiplier-convolution-commutation}
T_m^{\theta_3} \bigl(x\star_{\theta_{1},\theta_{2}}y\bigr)
= \bigl(T_m^{\theta_{1}}x\bigr)\star_{\theta_{1},\theta_{2}}y
= x\star_{\theta_{1},\theta_{2}}\bigl(T_m^{\theta_{2}}y\bigr).
\end{equation}
All these identities hold in $\mathcal{S}'(\mathbb{T}_{\theta_3}^{d})$ and hence in every function or operator space in which the displayed terms are defined.
\end{proposition}

\begin{proof}
The Fourier coefficients of a distribution on a quantum torus have at most polynomial growth, and conversely every polynomially growing sequence determines a unique distribution; see \cite[Sections~2.1--2.2]{XiongXuYin}. Therefore, $\widehat{x}(n)\widehat{y}(n)$ again has at most polynomial growth, which proves that \eqref{eq:distributional-convolution} defines a distribution. For $L_{1}$-inputs, the Fourier formula in Proposition~\ref{prop:basic-properties} shows that this extension agrees with the convolution defined as the preadjoint of the comultiplication.

For every $n\in\mathbb{Z}^{d}$, the $n$-th Fourier coefficient of the left-hand side of \eqref{eq:two-multiplier-convolution-identity} is
\begin{equation*}
m_{1}(n)m_{2}(n)\widehat{x}(n)\widehat{y}(n),
\end{equation*}
which is also the $n$-th Fourier coefficient of its right-hand side. Fourier uniqueness in $\mathcal{S}'$ proves \eqref{eq:two-multiplier-convolution-identity}. Taking successively $m_{2}=1$ and $m_{1}=1$ gives \eqref{eq:multiplier-convolution-commutation}.
\end{proof}

We next record the Besov module estimate needed below.

\begin{lem}[Besov module estimate]
\label{lem:Besov-module-estimates}
Let $a\in L_{1}(\mathbb{T}_{\theta}^{d}).$ Let $s\in\mathbb{R}$ and suppose that $g\in B_{\infty,\infty}^s(\mathbb{T}^{d}).$ We have $g\star_{0,\theta}a \in B_{\infty,\infty}^s(\mathbb{T}_{\theta}^{d})$ and
\begin{equation}\label{eq:Besov-module}
\|g\star_{0,\theta}a\|_{B_{\infty,\infty}^s}
\leq \|g\|_{B_{\infty,\infty}^s}\|a\|_{1}.
\end{equation}
Here the same Littlewood--Paley cutoff is used in the two Besov norms.
\end{lem}
\begin{proof} We now prove \eqref{eq:Besov-module}. Let $\varphi$ be the cutoff used in the Littlewood--Paley decomposition. By Proposition~\ref{prop:multiplier-compatibility},
\begin{equation}\label{eq:Besov-block-convolution-commutation}
\widetilde{\varphi}_{j}\ast \bigl(g\star_{0,\theta}a\bigr)
= \bigl(\widetilde{\varphi}_{j}\ast g\bigr) \star_{0,\theta}a.
\end{equation}
Proposition \ref{prop:L1-Linfty} now gives
\begin{equation*}
\left\|
\widetilde{\varphi}_{j}\ast
\bigl(g\star_{0,\theta}a\bigr)
\right\|_{\infty}
\leq
\|\widetilde{\varphi}_{j}\ast g\|_{\infty}\|a\|_{1}.
\end{equation*}
Moreover,
\begin{equation*}
\left|
\widehat{g\star_{0,\theta}a}(0)
\right|
= |\widehat g(0)|\,|\widehat a(0)|
\leq |\widehat g(0)|\|a\|_{1}.
\end{equation*}
Multiplying the dyadic estimate by $2^{js}$ and taking the supremum over $j\geq0$ proves \eqref{eq:Besov-module}.
\end{proof}

\begin{lem}\label{lem:Gk-Besov} Let $k\geq1$ be an integer. Set
\begin{equation*}
    G_k=\sum_{0\neq n\in\mathbb{Z}^{d}}\frac1{\sum_{j=1}^{d}n_{j}^{2k}}e_n.
\end{equation*} 
We have
\begin{equation*}
    G_k\in B_{\infty,\infty}^{2k-d}(\mathbb{T}^{d}),\ \  k\in\mathbb{N}.
\end{equation*}
\end{lem}
\begin{proof} Let $\varphi\in C_c^{\infty}(\mathbb{R}^{d})$ be the annular cutoff used in the Littlewood--Paley decomposition and define
\begin{equation*}
    \Delta_mG_k=\sum_{0\neq n\in\mathbb{Z}^{d}}\frac{\varphi(\frac{n}{2^m})}{\sum_{j=1}^{d}n_{j}^{2k}}e_n,\ \  m\in\mathbb{Z}_+.
\end{equation*}
We have
\begin{equation*}
    \begin{split}
        \|\Delta_mG_k\|_{\infty} & \leq\sum_{0\neq n\in\mathbb{Z}^{d}}\frac{|\varphi|(\frac{n}{2^m})}{\sum_{j=1}^{d}n_{j}^{2k}}\leq\|\varphi\|_{\infty}\sum_{\substack{n\in\mathbb{Z}^{d}\\ 2^{m-2}\leq |n|\leq 2^{m+3}}}\frac1{\sum_{j=1}^{d}n_{j}^{2k}}\\
        & \leq d^{k-1}\|\varphi\|_{\infty}\sum_{\substack{n\in\mathbb{Z}^{d}\\ 2^{m-2}\leq |n|\leq 2^{m+3}}}|n|^{-2k}.
    \end{split}
\end{equation*}
If $2k>d,$ then we estimate
\begin{equation*}
    \sum_{\substack{n\in\mathbb{Z}^{d}\\ 2^{m-2}\leq |n|\leq 2^{m+3}}}|n|^{-2k}\leq \sum_{\substack{n\in\mathbb{Z}^{d}\\ 2^{m-2}\leq |n|}}|n|^{-2k}\leq c_{d,k}2^{m(d-2k)}.
\end{equation*}
If $2k<d,$ then we estimate
\begin{equation*}
    \sum_{\substack{n\in\mathbb{Z}^{d}\\ 2^{m-2}\leq |n|\leq 2^{m+3}}}|n|^{-2k}\leq \sum_{\substack{0\neq n\in\mathbb{Z}^{d}\\ |n|\leq 2^{m+3}}}|n|^{-2k}\leq c_{d,k}2^{m(d-2k)}.
\end{equation*} 
If $2k=d,$ then we estimate 
\begin{equation*}
    \sum_{\substack{n\in\mathbb{Z}^{d}\\ 2^{m-2}\leq |n|\leq 2^{m+3}}}|n|^{-2k}=\sum_{\substack{n\in\mathbb{Z}^{d}\\ 2^{m-2}\leq |n|\leq 2^{m+3}}}|n|^{-d}\leq c_{d,k}.
\end{equation*}
Appealing to the definition of $B^{2k-d}_{\infty,\infty}(\mathbb{T}^{d})$, we complete the proof.
\end{proof}

We first construct the scalar kernels needed for the Sobolev--Besov embedding.

\begin{lem}[Green kernels for the pure-derivative reconstruction]
\label{lem:Besov-reconstruction-kernels}
Let $k\geq1$ be an integer. For $1\leq\ell\leq d,$ set
\begin{equation}
    G_{k,\ell}=\sum_{0\neq n\in\mathbb{Z}^{d}}\frac{n_{\ell}^{k}}{\sum_{j=1}^{d}n_{j}^{2k}}e_n.
\end{equation} 
We have
\begin{equation}
    G_{k,\ell} \in B_{\infty,\infty}^{k-d}(\mathbb{T}^{d}),\ \  1\leq\ell\leq d,\ \  k\in\mathbb{N}.
\end{equation}
\end{lem}\begin{proof}
By the definitions of $G_k$ and $G_{k,\ell}$, we have 
\begin{equation*}
G_{k,\ell}=D_{\ell}^{k}G_k
\end{equation*}
in $\mathcal{S}'(\mathbb{T}^{d})$. By
Lemma~\ref{lem:Gk-Besov},
\begin{equation*}
G_k\in B_{\infty,\infty}^{2k-d}(\mathbb{T}^{d}).
\end{equation*}
Let
\begin{equation*}
a=(0,\ldots,0,k,0,\ldots,0)\in\mathbb{N}_{0}^{d},
\end{equation*}
where the only nonzero entry is in the $\ell$-th position. Then
\begin{equation*}
D^{a}=D_{\ell}^{k} \ \text{ and }  |a|_{1}=k.
\end{equation*}
The differentiation property for Besov spaces \cite[Theorem~3.7(ii)]{XiongXuYin}, applied with $\theta=0$, $p=q=\infty$, $s=2k-d$ and the preceding multi-index $a$, gives
\begin{equation*}
D_{\ell}^{k}: B_{\infty,\infty}^{2k-d}(\mathbb{T}^{d})
\to 
B_{\infty,\infty}^{k-d}(\mathbb{T}^{d}).
\end{equation*}
Consequently,
\begin{equation*}
G_{k,\ell} = D_{\ell}^{k}G_k \in B_{\infty,\infty}^{k-d}(\mathbb{T}^{d}),
\end{equation*}
as required.
\end{proof}

The preceding lemma gives a scalar resolution of the identity adapted to the pure derivatives of order $k$.

\begin{thm}[Endpoint Sobolev--Besov embedding]
\label{thm:W1-Besov-endpoint}
Let $d\geq1$, let $\theta$ be a real skew-symmetric $d\times d$ matrix, and let $k\geq1$ be an integer. Then
\begin{equation}\label{eq:W1-Besov-endpoint}
W_{1}^{k}(\mathbb{T}_{\theta}^{d})
\hookrightarrow
B_{\infty,\infty}^{k-d}(\mathbb{T}_{\theta}^{d}).
\end{equation}
More precisely, there exists a constant $C_{d,k}>0,$ independent of $\theta$, such that
\begin{equation}\label{eq:W1-Besov-endpoint-estimate}
\|x\|_{B_{\infty,\infty}^{k-d}(\mathbb{T}_{\theta}^{d})}
\leq C_{d,k} \left( \|x\|_{1} + \sum_{\ell=1}^{d}\|\partial_{\ell}^{k}x\|_{1} \right)
\end{equation}
for every $x\in W_{1}^{k}(\mathbb{T}_{\theta}^{d})$.
\end{thm}
\begin{proof} We have
$$\sum_{\ell=1}^{d}\widehat{G_{k,\ell}}(n)\cdot n_{\ell}^{k}=1,\ \  0\neq n\in\mathbb{Z}^{d}.$$
Therefore,
$$x=\widehat{x}(0)\mathbf{1}+\sum_{\ell=1}^{d}G_{k,\ell}\star_{0,\theta}D_{\ell}^{k}x.$$
Consequently,
\begin{equation*}
    \begin{split}
        \|x\|_{B_{\infty,\infty}^{k-d}(\mathbb{T}_{\theta}^{d})} & \leq |\hat{x}(0)|+\sum_{\ell=1}^{d}\|G_{k,\ell}\star_{0,\theta}D_{\ell}^{k}x\|_{B_{\infty,\infty}^{k-d}(\mathbb{T}_{\theta}^{d})}\\
        & \leq|\hat{x}(0)|+\sum_{\ell=1}^{d}\|G_{k,\ell}\|_{B_{\infty,\infty}^{k-d}(\mathbb{T}^{d})}\|D_{\ell}^{k}x\|_{L_{1}(\mathbb{T}^{d}_{\theta})}.
    \end{split}
\end{equation*} 
Since $|\widehat{x}(0)|\leq\|x\|_{1}$ and $\|D_{\ell}^{k}x\|_{1}=\|\partial_{\ell}^{k}x\|_{1}$, the assertion follows from Lemma~\ref{lem:Besov-reconstruction-kernels}.
\end{proof}

\begin{lem}[$L_{1}$-Besov regularity of the reconstruction kernels]
\label{lem:Gkl-B1-infinity}
Let $k\geq1$ and $1\leq\ell\leq d$. Then
\begin{equation*}
G_{k,\ell}\in B_{1,\infty}^{k}(\mathbb{T}^{d}).
\end{equation*}
More precisely,
\begin{equation}\label{eq:Gkl-dyadic-L1}
\left\| \widetilde{\varphi}_{j}\ast G_{k,\ell} \right\|_{L_{1}(\mathbb{T}^{d})}
\leq C_{d,k,\varphi}2^{-jk}, \ \  j\geq0.
\end{equation}
\end{lem}

\begin{proof}
Define
\begin{equation*}
\psi_{k,\ell}(\xi)
=
\begin{cases}
\displaystyle
\frac{\xi_{\ell}^{k}\varphi(\xi)}
{\sum_{m=1}^{d}\xi_m^{2k}},
& \xi\neq0, \\
0,
& \xi=0.
\end{cases}
\end{equation*}
Since $\varphi$ is smooth and supported in an annulus separated from
the origin,
\begin{equation*}
\psi_{k,\ell}\in C_c^{\infty}(\mathbb{R}^{d}).
\end{equation*}
For every $n\in\mathbb{Z}^{d}$,
\begin{align*}
\widehat{\widetilde{\varphi}_{j}\ast G_{k,\ell}}(n)
= \varphi(2^{-j}n) \frac{n_{\ell}^{k}}{\sum_{m=1}^{d}n_m^{2k}} = 2^{-jk}\psi_{k,\ell}(2^{-j}n).
\end{align*}

Put
\begin{equation*}
K_{j,\ell}(t) = \sum_{n\in\mathbb{Z}^{d}} \psi_{k,\ell}(2^{-j}n)e^{in\cdot t}.
\end{equation*}
By the Poisson summation formula, $K_{j,\ell}$ is, up to a fixed normalization constant, the periodization of
\begin{equation*}
t\mapsto 2^{jd}\check{\psi}_{k,\ell}(2^jt).
\end{equation*}
Consequently, if $Q=[0,2\pi]^{d}$, then
\begin{align*}
\|K_{j,\ell}\|_{L_{1}(\mathbb{T}^{d})}
&\leq
C_d\sum_{m\in\mathbb{Z}^{d}}
\int\limits_Q
2^{jd}
\left|
\check{\psi}_{k,\ell}
\bigl(2^j(t+2\pi m)\bigr)
\right|\,dt \\
&= C_d \|\check{\psi}_{k,\ell}\|_{L_{1}(\mathbb{R}^{d})} \leq C_{d,k,\varphi}.
\end{align*}
Since
\begin{equation*}
\widetilde{\varphi}_{j}\ast G_{k,\ell} = 2^{-jk}K_{j,\ell},
\end{equation*}
we obtain \eqref{eq:Gkl-dyadic-L1}. The low-frequency block of $G_{k,\ell}$ is a trigonometric polynomial and hence belongs to $L_{1}(\mathbb{T}^{d})$. Therefore, $G_{k,\ell}\in B_{1,\infty}^{k}(\mathbb{T}^{d}).$
\end{proof}

\begin{thm}[Pure-derivative Besov embedding at general integrability]
\label{thm:Wp-Besov-odd-endpoint}
Let $d,k\geq1$ be integers and let $1\leq p\leq\infty$. Let $p'$ be
the conjugate exponent of $p$. Then
\begin{equation}\label{eq:Wp-Besov-odd-endpoint}
W_{p}^{k}(\mathbb{T}_{\theta}^{d})
\hookrightarrow
B_{\infty,\infty}^{k-\frac{d}{p}}(\mathbb{T}_{\theta}^{d}).
\end{equation}
More precisely, with the same Littlewood--Paley cutoff on the
classical and quantum tori,
\begin{equation}\label{eq:Wp-Besov-odd-bound}
\|x\|_{B_{\infty,\infty}^{k-\frac{d}{p}}}
\leq
C_{d,k}
\left(
\|x\|_{p}+
\sum_{\ell=1}^{d}\|\partial_\ell^{k}x\|_{p}
\right).
\end{equation}
The constant is independent of $\theta$. In particular, this covers
the case of odd $k$ in \cite[Remark~6.8]{XiongXuYin}.
\end{thm}

\begin{proof}
The reconstruction formula gives
\begin{equation*}
x
=
\widehat{x}(0)\mathbf{1}
+
\sum_{\ell=1}^{d}
G_{k,\ell}\star_{0,\theta}D_{\ell}^{k}x.
\end{equation*}
Hence, by Proposition~\ref{prop:associative-commutative},
\begin{equation*}
\widetilde{\varphi}_{j}\star_{0,\theta}x
=
\sum_{\ell=1}^{d}
\bigl(\widetilde{\varphi}_{j}\ast G_{k,\ell}\bigr)
\star_{0,\theta}D_{\ell}^{k}x.
\end{equation*}
Applying Theorem~\ref{thm:Young} with $r=\infty$, we obtain
\begin{align*}
\left\|
\widetilde{\varphi}_{j}\star_{0,\theta}x
\right\|_{\infty}
\leq
\sum_{\ell=1}^{d}
\left\|
\widetilde{\varphi}_{j}\ast G_{k,\ell}
\right\|_{p'}
\left\|D_{\ell}^{k}x\right\|_{p}.
\end{align*}

By Lemma~\ref{lem:Besov-reconstruction-kernels},
\begin{equation*}
\left\|
\widetilde{\varphi}_{j}\ast G_{k,\ell}
\right\|_{\infty}
\leq
C_{d,k,\varphi}2^{-j(k-d)},
\end{equation*}
whereas Lemma~\ref{lem:Gkl-B1-infinity} gives
\begin{equation*}
\left\|
\widetilde{\varphi}_{j}\ast G_{k,\ell}
\right\|_{1}
\leq
C_{d,k,\varphi}2^{-jk}.
\end{equation*}
Interpolating these two estimates, with the usual endpoint
conventions, yields
\begin{align*}
\left\|
\widetilde{\varphi}_{j}\ast G_{k,\ell}
\right\|_{p'}
&\leq
\left\|
\widetilde{\varphi}_{j}\ast G_{k,\ell}
\right\|_{1}^{\frac1{p'}}
\left\|
\widetilde{\varphi}_{j}\ast G_{k,\ell}
\right\|_\infty^{1-\frac1{p'}} \\
&\leq
C_{d,k,\varphi}
2^{-\frac{jk}{p'}}
2^{-j(k-d)(1-\frac1{p'})} \\
&=
C_{d,k,\varphi}
2^{-j(k-\frac{d}{p})}.
\end{align*}
Therefore,
\begin{equation*}
\left\|
\widetilde{\varphi}_{j}\star_{0,\theta}x
\right\|_{\infty}
\leq
C_{d,k,\varphi}
2^{-j(k-\frac{d}{p})}
\sum_{\ell=1}^{d}
\left\|D_{\ell}^{k}x\right\|_{p}.
\end{equation*}
Finally, $|\widehat{x}(0)|\leq\|x\|_{p}$ and $\|D_{\ell}^{k}x\|_{p} = \|\partial_{\ell}^{k}x\|_{p}.$ The definition of $B_{\infty,\infty}^{k-\frac{d}{p}}(\mathbb{T}_{\theta}^{d})$ now gives \eqref{eq:Wp-Besov-odd-bound}.
\end{proof}

We now turn to the critical $L_{\infty}$-embedding. Here the appropriate resolution of the identity is indexed by the nonempty subsets of $\{1,\ldots,d\}$ and involves mixed derivatives.

\begin{lem}\label{lem:sawtooth-reconstruction-kernels} Set
$$h=\sum_{\substack{n\in\mathbb{Z}^{d}\\ n_{j}\neq0 \textnormal{ for } 1\leq j\leq d}}\prod_{1\leq j\leq d}\frac1{n_{j}}\cdot e_n.$$
We have
$$\|h\|_{\infty}=\pi^{d}.$$
\end{lem}
\begin{proof} Denote
$$h_0=\sum_{0\neq m\in\mathbb{Z}}\frac1me_m.$$
We have $h=h_0^{\otimes d}$ and, therefore,
$$\|h\|_{\infty}=\|h_0\|_{\infty}^{d}.$$
It remains to note that $h_0$ is the $2\pi$-periodic extension of the function $t\to i(\pi-t),$ $t\in(0,2\pi).$ Indeed, with normalized Haar measure,
\begin{equation*}
\widehat{i(\pi-t)}(m)=\frac{1}{m},
\ \  m\neq0,
\end{equation*}
and its zeroth Fourier coefficient vanishes. In particular,
$\|h_0\|_{\infty}=\pi.$
\end{proof}

\begin{thm}[Critical endpoint Sobolev embedding]
\label{thm:critical-W1-Linfty}
Let $d\geq1$ and let $\theta$ be a real skew-symmetric $d\times d$ matrix. We have
\begin{equation}\label{eq:critical-W1-Linfty}
W_{1}^{d}(\mathbb{T}_{\theta}^{d})
\hookrightarrow
L_{\infty}(\mathbb{T}_{\theta}^{d}).
\end{equation}
\end{thm}
\begin{proof} Let $h$ be as in Lemma \ref{lem:sawtooth-reconstruction-kernels}. We have
$$\hat{h}(n)\cdot \prod_{j=1}^{d}n_{j}=
\begin{cases}
1,& n_{j}\neq0\textnormal{ for }1\leq j\leq d\\
0,&\mbox{otherwise}
\end{cases}
$$
Therefore,
$$\sum_{\substack{n\in\mathbb{Z}^{d}\\ n_{j}\neq0\textnormal{ for }1\leq j\leq d}}\hat{x}(n)e_n^{\theta}=h\star_{0,\theta}(D_{1}\cdots D_dx).$$
It follows that
\begin{equation*}
    \begin{split}
        \Big\|\sum_{\substack{n\in\mathbb{Z}^{d}\\ n_{j}\neq0\textnormal{ for }1  \leq j\leq d}}\hat{x}(n)e_n^{\theta}\Big\|_{L_{\infty}(\mathbb{T}^{d}_{\theta})} & \leq\|h\|_{L_{\infty}(\mathbb{T}^{d})}\|D_{1}\cdots D_dx\|_{L_{1}(\mathbb{T}^{d}_{\theta})}\\
        & =\pi^{d}\|D_{1}\cdots D_dx\|_{L_{1}(\mathbb{T}^{d}_{\theta})}.
    \end{split}
\end{equation*}
It remains to note that
\begin{equation*}
    x=\sum_{S\subset\{1,\cdots,d\}}\sum_{\substack{n\in\mathbb{Z}^{d}\\ n_{j}\neq0\textnormal{ iff }j\in S}}\hat{x}(n)e_n^{\theta}.
\end{equation*} 
Denote the inner sum by $x_S$. The summand corresponding to $S=\{1,\cdots,d\}$ is treated above, while $x_{\varnothing}=\widehat{x}(0)\mathbf{1}$. If $\varnothing\neq S\subsetneq\{1,\cdots,d\}$, then $x_S\in W_{1}^{d}(\mathbb{T}^{|S|}_{\theta_S})$, where $\theta_S=(\theta_{k,l})_{k,l\in S}$. Indeed, the Fourier projection onto the frequencies with $\operatorname{supp}_0(n)=S$ is a finite alternating sum of the trace-preserving conditional expectations onto coordinate subtori; hence it is bounded on $L_{1}$ and commutes with the canonical derivations.

Applying Theorem~\ref{thm:W1-Besov-endpoint} in dimension $|S|$ and with the order $k=d$, we obtain
\begin{equation*}
x_S\in B_{\infty,\infty}^{d-|S|} (\mathbb{T}^{|S|}_{\theta_S}).
\end{equation*}
Since $d-|S|>0$, the Littlewood--Paley decomposition and a geometric series give
\begin{align*}
\|x_S\|_{\infty} \leq |\widehat{x_S}(0)|+
\sum_{j\geq0}\|\widetilde\varphi_{j}\ast x_S\|_{\infty} \leq
\left(1+\sum_{j\geq0}2^{-j(d-|S|)}\right)
\|x_S\|_{B_{\infty,\infty}^{d-|S|}}.
\end{align*}
The coordinate projections are bounded on $L_{1}$, and each lower-dimensional Besov estimate is bounded by the $W_{1}^{d}$ norm of $x$. Since the sum over $S$ is finite, the preceding estimates give $\|x\|_\infty\leq C_d\|x\|_{W_{1}^{d}}$, which proves \eqref{eq:critical-W1-Linfty}.
\end{proof}

\begin{rem} The embedding $W_{1}^{k} (\mathbb{T}^{d}_{\theta}) \hookrightarrow W_\infty^{k-d} (\mathbb{T}^{d}_{\theta})$ for integers $k\geq d$ can be proved similarly.
\end{rem}

\appendix
\section{Construction of the parameter-reversing anti-isomorphism}\label{app:parameter-reversal}
We give the steps used in Lemma~\ref{lem:parameter-reversal}. Let $\mathcal{P}_{\theta}$ denote the algebra of Fourier polynomials on $\mathbb{T}_{\theta}^{d}$.

\begin{lem}[Algebraic step]\label{lem:kappa-polynomials}
The linear map
\begin{equation}\label{eq:kappa-polynomial-definition}
\kappa_{\theta}\left(\sum_{n}a_{n}e_{n}^{\theta}\right) = \sum_{n}a_{n}e_{-n}^{-\theta}
\end{equation}
is a unital bijective $*$-anti-homomorphism from $\mathcal{P}_{\theta}$ onto $\mathcal{P}_{-\theta}$.  Its inverse is $\kappa_{-\theta}$.
\end{lem}

\begin{proof}
Skew-symmetry of $\theta$ gives
\begin{align*}
\sigma_{-\theta}(-m,-n)
&= \exp\bigl(\pi i\langle-m,(-\theta)(-n)\rangle\bigr) \\
&= \exp\bigl(-\pi i\langle m,\theta n\rangle\bigr)
= \exp\bigl(\pi i\langle n,\theta m\rangle\bigr) \\
&= \sigma_{\theta}(n,m).
\end{align*}
Consequently,
\begin{align*}
\kappa_{\theta}(e_{n}^{\theta}e_{m}^{\theta})
&= \sigma_{\theta}(n,m)e_{-(n+m)}^{-\theta} = \sigma_{-\theta}(-m,-n)e_{-(m+n)}^{-\theta} \\
&= e_{-m}^{-\theta}e_{-n}^{-\theta}
= \kappa_{\theta}(e_{m}^{\theta}) \kappa_{\theta}(e_{n}^{\theta}).
\end{align*}
Also,
\begin{equation*}
\kappa_{\theta}\bigl((e_{n}^{\theta})^{*}\bigr)
= e_{n}^{-\theta}
= \kappa_{\theta}(e_{n}^{\theta})^{*}.
\end{equation*}
Linearity proves anti-multiplicativity.  It also proves preservation of the involution on Fourier polynomials.  Finally, $\kappa_{-\theta}\kappa_{\theta}$ fixes every Fourier generator, and the same argument with $\theta$ and $-\theta$ interchanged gives the reverse composition.
\end{proof}

\begin{lem}[$C^{*}$-algebraic step]\label{lem:kappa-Cstar}
The map in Lemma~\ref{lem:kappa-polynomials} extends uniquely to an isometric $*$-anti-isomorphism
\begin{equation*}
\kappa_{\theta}: C(\mathbb{T}_{\theta}^{d})
\to  C(\mathbb{T}_{-\theta}^{d}).
\end{equation*}
It preserves the canonical traces:
\begin{equation}\label{eq:kappa-trace-Cstar}
\tau_{-\theta}(\kappa_{\theta}(a))=\tau_{\theta}(a), \ \  a\in C(\mathbb{T}_{\theta}^{d}).
\end{equation}
\end{lem}

\begin{proof}
Let $C(\mathbb{T}_{-\theta}^{d})^{\mathrm{op}}$ be the opposite $C^{*}$-algebra.  It has the same vector space, norm and involution as $C(\mathbb{T}_{-\theta}^{d})$, while its multiplication is
\begin{equation*}
a\mathbin{\cdot_{\mathrm{op}}}b=ba.
\end{equation*}
Put $v_{n}=e_{-n}^{-\theta}$ and regard $v_{n}$ as an element of the opposite algebra.  Lemma~\ref{lem:kappa-polynomials} gives
\begin{equation*}
v_{n}\mathbin{\cdot_{\mathrm{op}}}v_{m} = \sigma_{\theta}(n,m)v_{n+m}, \ \  v_{n}^{*}=v_{-n}.
\end{equation*}
Thus $(v_{n})$ satisfies the defining twisted relations.  For the conventions, see \cite[Definition~3.13]{LeeXiong} and \cite[Section~1.2]{XiongXuYin}.  The universal property gives a unital $*$-homomorphism
\begin{equation*}
C(\mathbb{T}_{\theta}^{d})
\to  C(\mathbb{T}_{-\theta}^{d})^{\mathrm{op}}.
\end{equation*}
Interpreted with the original multiplication in the target, this map is contractive and reverses products:
\begin{equation*}
\kappa_{\theta}: C(\mathbb{T}_{\theta}^{d}) \to  C(\mathbb{T}_{-\theta}^{d}).
\end{equation*}
The construction with $-\theta$ gives $\kappa_{-\theta}$.  The two compositions are the identity on the norm-dense Fourier-polynomial algebras, hence on the full $C^{*}$-algebras.  Both maps are contractive, so they are isometric.

For every $n\in\mathbb{Z}^{d}$,
\begin{equation*}
\tau_{-\theta}(\kappa_{\theta}(e_{n}^{\theta}))
= \tau_{-\theta}(e_{-n}^{-\theta})
= \delta_{n,0}
= \tau_{\theta}(e_{n}^{\theta}).
\end{equation*}
Linearity, norm continuity and density of $\mathcal{P}_{\theta}$ prove \eqref{eq:kappa-trace-Cstar}.
\end{proof}

We isolate the passage from $C^{*}$-algebras to finite von Neumann algebras.

\begin{lem}[Normal extension]\label{lem:trace-preserving-normal-extension}
Let $(M,\tau_{M})$ and $(N,\tau_{N})$ be finite von Neumann algebras with faithful normal tracial states.  Let $A\subset M$ and $B\subset N$ be unital weak-$*$ dense $C^{*}$-subalgebras.  Every trace-preserving $*$-isomorphism $\alpha:A\to B$ has a unique trace-preserving normal $*$-isomorphism extension $\overline\alpha:M\to N$.
\end{lem}

\begin{proof}
Trace preservation makes $a\mapsto\alpha(a)$ an isometry from $A\subset L_{2}(M,\tau_{M})$ onto $B\subset L_{2}(N,\tau_{N})$.  Kaplansky's density theorem \cite[Theorem~1]{Kaplansky1951}, see also \cite[Theorem~II.4.8]{Takesaki}, implies that $A$ and $B$ are dense in the respective $L_{2}$-spaces.  Hence the isometry extends to a unitary
\begin{equation*}
U:L_{2}(M,\tau_{M})\to  L_{2}(N,\tau_{N}).
\end{equation*}
Let $\lambda_{M}$ and $\lambda_{N}$ be the left multiplication representations.  For $a,c\in A$, 
\begin{equation*}
U\lambda_{M}(a)c = U(ac) =\alpha(a)\alpha(c) =\lambda_{N}(\alpha(a))U(c).
\end{equation*}
Thus
\begin{equation*}
U\lambda_{M}(a)U^{*}=\lambda_{N}(\alpha(a)), \ \  a\in A.
\end{equation*}
Taking von Neumann closures gives $U\lambda_{M}(M)U^{*}=\lambda_{N}(N)$.  Therefore
\begin{equation*}
\lambda_{N}(\overline\alpha(x)) = U\lambda_{M}(x)U^{*},
\ \  x\in M,
\end{equation*}
defines a normal $*$-isomorphism extending $\alpha$.  Since $U\1=\1$, it preserves the trace.  Uniqueness follows from normality and the weak-$*$ density of $A$.
\end{proof}

\begin{lem}[Von Neumann algebraic step]\label{lem:kappa-vna}
The map in Lemma~\ref{lem:kappa-Cstar} extends uniquely to a trace-preserving normal linear $*$-anti-isomorphism
\begin{equation*}
\kappa_{\theta}: L_{\infty}(\mathbb{T}_{\theta}^{d}) \to  L_{\infty}(\mathbb{T}_{-\theta}^{d}).
\end{equation*}
Its inverse is $\kappa_{-\theta}$, and
\begin{equation}\label{eq:kappa-trace}
\tau_{-\theta}(\kappa_{\theta}(a))=\tau_{\theta}(a), \ \  a\in L_{\infty}(\mathbb{T}_{\theta}^{d}).
\end{equation}
\end{lem}

\begin{proof}
Regard $\kappa_{\theta}$ as a trace-preserving $*$-isomorphism
\begin{equation*}
C(\mathbb{T}_{\theta}^{d}) \to  C(\mathbb{T}_{-\theta}^{d})^{\mathrm{op}}.
\end{equation*}
Apply Lemma~\ref{lem:trace-preserving-normal-extension} to the canonical tracial von Neumann closures.  Returning to the original multiplication in the target gives the required normal $*$-anti-isomorphism.  The inverse relation and trace identity hold on the weak-$*$ dense $C^{*}$-subalgebras and therefore pass to the normal extensions.
\end{proof}

\begin{lem}[$L_{p}$ and measurable-operator step]\label{lem:kappa-Lp}
For every $1\leq p<\infty$, $\kappa_{\theta}$ extends uniquely to a surjective linear isometry
\begin{equation*}
\kappa_{\theta}: L_{p}(\mathbb{T}_{\theta}^{d}) \to  L_{p}(\mathbb{T}_{-\theta}^{d}).
\end{equation*}
For $x\in L_{p}(\mathbb{T}_{\theta}^{d})$ and $a\in L_{\infty}(\mathbb{T}_{\theta}^{d})$,
\begin{equation}\label{eq:kappa-module-identities}
\kappa_{\theta}(ax)=\kappa_{\theta}(x)\kappa_{\theta}(a),
\ \  \kappa_{\theta}(xa)=\kappa_{\theta}(a)\kappa_{\theta}(x).
\end{equation}
When $p=1$,
\begin{equation}\label{eq:kappa-trace-L1}
\tau_{-\theta}(\kappa_{\theta}(x))=\tau_{\theta}(x).
\end{equation}
Moreover, $\kappa_{\theta}$ extends uniquely to
\begin{equation*}
\kappa_{\theta}: S(\mathbb{T}_{\theta}^{d}) \to  S(\mathbb{T}_{-\theta}^{d})
\end{equation*}
and preserves generalized singular values, as in \eqref{eq:kappa-mu}.
\end{lem}

\begin{proof}
For $a\in L_{\infty}(\mathbb{T}_{\theta}^{d})$,
\begin{align*}
|\kappa_{\theta}(a)|^{2} =\kappa_{\theta}(a)^{*}\kappa_{\theta}(a)
= \kappa_{\theta}(a^{*})\kappa_{\theta}(a)
= \kappa_{\theta}(aa^{*}).
\end{align*}
The von Neumann algebra generated by $aa^{*}$ is commutative, so order reversal has no effect there.  Functional calculus yields
\begin{equation}\label{eq:kappa-absolute-value}
|\kappa_{\theta}(a)|^{p} = \kappa_{\theta}(|a^{*}|^{p}), \ \  1\leq p<\infty.
\end{equation}
Using \eqref{eq:kappa-trace} and traciality,
\begin{align*}
\|\kappa_{\theta}(a)\|_{p}^{p}
= \tau_{-\theta}\bigl(|\kappa_{\theta}(a)|^{p}\bigr)
= \tau_{\theta}(|a^{*}|^{p})
= \tau_{\theta}(|a|^{p})
= \|a\|_{p}^{p}.
\end{align*}
The penultimate equality follows from the polar decomposition and traciality.  If $x=u|x|\in L_{p}$, then
\begin{equation*}
x_{N}=u|x|\1_{[0,N]}(|x|)
\end{equation*}
is bounded and
\begin{equation*}
\|x-x_{N}\|_{p}^{p} = \tau_{\theta}\bigl(|x|^{p}\1_{(N,\infty)}(|x|)\bigr) \to 0.
\end{equation*}

Thus bounded operators are dense in $L_{p}$, and $\kappa_{\theta}$ extends isometrically.  The inverse relation gives surjectivity.  Approximating $x$ by bounded operators also proves \eqref{eq:kappa-module-identities} and, for $p=1$, \eqref{eq:kappa-trace-L1}.

View the normal anti-isomorphism as an isomorphism into the opposite algebra.  Propositions~2.9.3 and 3.3.10(i) of \cite{DoddsDePagterSukochev} give its unique extension to measurable operators.  Proposition~3.3.10(ii) of the same reference gives preservation of singular values, proving \eqref{eq:kappa-mu}.
\end{proof}

\begin{lem}[Fourier step]\label{lem:kappa-Fourier-step}
For every $x\in L_{1}(\mathbb{T}_{\theta}^{d})$ and $n\in\mathbb{Z}^{d}$,
\begin{equation*}
\widehat{\kappa_{\theta}(x)}(n)=\widehat{x}(-n).
\end{equation*}
\end{lem}

\begin{proof}
Equations \eqref{eq:kappa-basis}, \eqref{eq:kappa-module-identities} and \eqref{eq:kappa-trace-L1} give 
\begin{align*}
\widehat{\kappa_{\theta}(x)}(n)
&= \tau_{-\theta}\bigl(\kappa_{\theta}(x)e_{-n}^{-\theta}\bigr)
= \tau_{-\theta} \bigl(\kappa_{\theta}(e_{n}^{\theta}x)\bigr)
= \tau_{\theta}(e_{n}^{\theta}x) \\
&= \tau_{\theta}(xe_{n}^{\theta}) = \widehat{x}(-n).
\end{align*}
The cyclicity used here is valid for the product of an $L_{\infty}$-operator and an $L_{1}$-operator.  It follows by approximating $x$ in $L_{1}$ by bounded operators.
\end{proof}

\begin{proof}[Proof of Lemma~\ref{lem:parameter-reversal}]
The successive assertions follow from Lemmas~\ref{lem:kappa-polynomials}, \ref{lem:kappa-Cstar}, \ref{lem:kappa-vna}, \ref{lem:kappa-Lp} and \ref{lem:kappa-Fourier-step}.
\end{proof}

{\bf Conflict of Interest Statement:} The author declares that there is no conflict of interest.

{\bf Data Availability Statement:} No datasets were generated or analyzed during the current study.

{\bf Acknowledgement}

%The work was partially supported by the grant No. AP23487088 of the Science Committee of the Ministry of Science and Higher
%Education of the Republic of Kazakhstan.
The authors would like to thank Professor Xiao Xiong  for useful discussions and for confirming the answer to the question posed in their joint paper with Professors Xu and Yin in \cite{XiongXuYin}. The authors would also like to express their special thanks to Dr. Edward McDonald for his helpful discussions and valuable comments.
The second and third authors would also like to thank Professor Jinghao Huang  for inviting them to the Institute for Advanced Study in Mathematics (IASM) at Harbin Institute of Technology and for providing a warm and stimulating atmosphere in which the basic ideas of this paper emerged.
The work was partially supported by the Australian Research Council. F.S. was supported by ARC grant FL17010005 and ARC grant DP230100434. K.T. was partially supported by ARC grant FL17010005. No AI tools were used to develop the arguments, mathematical content, or citations of this article. The free version of ChatGPT was used solely for language checking, including grammar, spelling, and wording.
%D.Z. was supported by ARC grant DP230100434. 

%Authors thank the anonymous referees for reading the paper and providing thoughtful comments, which improved the exposition of the paper.

\end{document}